\documentclass[11pt,letterpaper,reqno]{amsart}

\usepackage{amssymb}
\usepackage{amsmath}
\usepackage{amsthm}
\usepackage{amsfonts}
\usepackage{bbm}
\usepackage{esint}

\usepackage{bookmark}
\usepackage{hyperref}
\hypersetup{pdfstartview={FitH}}

\newtheorem{theorem}{Theorem}

\numberwithin{equation}{section}

\begin{document}

\title[Density theorems for anisotropic point configurations]{Density theorems for anisotropic point configurations}

\author[Vjekoslav Kova\v{c}]{Vjekoslav Kova\v{c}}
\address{Vjekoslav Kova\v{c}, Department of Mathematics, Faculty of Science, University of Zagreb, Bijeni\v{c}ka cesta 30, 10000 Zagreb, Croatia}
\email{vjekovac@math.hr}




\begin{abstract}
Several results in the existing literature establish Euclidean density theorems of the following strong type. These results claim that every set of positive upper Banach density in the Euclidean space of an appropriate dimension contains isometric copies of all sufficiently large elements of a prescribed family of finite point configurations. So far, all results of this type discussed linear isotropic dilates of a fixed point configuration. In this paper we initiate the study of analogous density theorems for families of point configurations generated by anisotropic dilations, i.e., families with power-type dependence on a single parameter interpreted as their size. More specifically, here we prove nonisotropic power-type generalizations of a result by Bourgain on vertices of a simplex, a result by Lyall and Magyar on vertices of a rectangular box, and a result on distance trees, which is a particular case of the treatise of distance graphs by Lyall and Magyar. Another source of motivation for this paper is providing additional evidence for the versatility of the approach stemming from the work of Cook, Magyar, and Pramanik and its modification used recently by Durcik and the present author. Finally, yet another purpose of this paper is to single out anisotropic multilinear singular integral operators associated with the above combinatorial problems, as they are interesting on their own.
\end{abstract}

\maketitle


\section{Introduction}
\subsection{Overview of previous results}
The topic of our interest are \emph{density theorems} within the subfield of the Euclidean Ramsey theory. Such theorems typically attempt to identify ``many'' dilates of a given finite point configuration in all sufficiently large measurable sets $A\subseteq\mathbb{R}^d$. Here the concept of largeness has to be interpreted in an appropriate measure-theoretic sense: by requiring that a certain density of $A$ is strictly positive. A general and very convenient notion of density is the \emph{upper Banach density}, defined for a measurable set $A\subseteq\mathbb{R}^d$ as
\[ \overline{\delta}(A) := \limsup_{R\rightarrow\infty}\sup_{x\in\mathbb{R}^d}\frac{|A\cap(x+[0,R]^d)|}{R^d}. \]

A quite strong density theorem for the simplest possible point configuration, namely the set of two points, was shown independently by Furstenberg, Katznelson, and Weiss \cite{FKW90:dist}, and Falconer and Marstrand \cite{FM86:dist}:
\begin{itemize}
\item[] \emph{For every measurable set $A\subseteq\mathbb{R}^2$ satisfying $\overline{\delta}(A)>0$ there is a positive number $\lambda_0=\lambda_0(A)$ such that for each $\lambda\in[\lambda_0,\infty)$ there exist points $x,x'\in A$ satisfying $|x-x'|=\lambda$.}
\end{itemize}
Here and in what follows, $|v|$ denotes the Euclidean norm of a vector $v\in\mathbb{R}^d$. The claim extends to higher dimensions, but we always consider only the smallest dimension in which the result is known to hold.
Bourgain \cite{B86:roth} generalized the above result to the set of vertices $\Delta\subseteq\mathbb{R}^{n}$ of a non-degenerate $n$-dimensional (i.e.\@ $(n+1)$-point) simplex:
\begin{itemize}
\item[] \emph{For every measurable set $A\subseteq\mathbb{R}^{n+1}$ satisfying $\overline{\delta}(A)>0$ there is a positive number $\lambda_0=\lambda_0(A,\Delta)$ such that for each $\lambda\in[\lambda_0,\infty)$ the set $A$ contains an isometric copy of $\lambda\Delta$.}
\end{itemize}
Note that the configuration $\Delta$ is initially given in the $n$-dimensional Euclidean space, while the ambient Euclidean space (i.e., the one containing the set $A$) has dimension $n+1$.
This dimensional increase is used in all known proofs of the aforementioned result, giving an additional ``degree of freedom,'' but at the time of writing it is not known if it is really necessary when $n\geq 2$.

The paper \cite{B86:roth} was very influential and it motivated a series of papers handling more complicated point configurations. Pursuing one possible direction, Lyall and Magyar initiated the study of density theorems for product-type point configurations. In \cite{LM16:prod} they considered Cartesian products $\Delta_1\times\Delta_2$, where both $\Delta_1$ and $\Delta_2$ are sets of vertices of non-degenerate simplices, while in \cite{LM19:hypergraphs} they extended their study to Cartesian products of finitely many such sets. An interesting (and already nontrivial) particular case is the set of vertices of an $n$-dimensional rectangular box,
$\Box = \{0,b_1\} \times \cdots \times \{0,b_n\} \subseteq \mathbb{R}^n$ for some $b_1,\ldots,b_n>0$.
One of the results by Lyall and Magyar, \cite[Theorem~1.1\,(i)]{LM19:hypergraphs}, reads as follows:
\begin{itemize}
\item[] \emph{For every measurable set $A\subseteq\mathbb{R}^{2}\times\cdots\times\mathbb{R}^{2}=(\mathbb{R}^2)^n$ satisfying $\overline{\delta}(A)>0$ there is a positive number $\lambda_0=\lambda_0(A,\Box)$ such that for each $\lambda\in[\lambda_0,\infty)$ the set $A$ contains an isometric copy of $\lambda\Box$ with sides parallel to the distinguished $2$-dimensional coordinate planes.
In other words, for each $\lambda\in[\lambda_0,\infty)$ one can find
$x_1,\ldots,x_n, y_1,\ldots,y_n \in\mathbb{R}^2$ satisfying
\[ \big\{ (x_1 + r_1 y_1, x_2 + r_2 y_2, \ldots, x_n + r_n y_n) : (r_1,\ldots,r_n)\in\{0,1\}^n \big\} \subseteq A \]
and $|y_k| = \lambda b_k$ for $k=1,2,\ldots,n$.}
\end{itemize}
In fact, Durcik and the present author first established a weaker result \cite[Theorem~1]{DK18:products}, with $(\mathbb{R}^2)^n$ replaced by $(\mathbb{R}^5)^n$, and then also reproved the above result \cite[Theorem~3]{DK20:Szemeredi}. The main concern of the paper \cite{DK20:Szemeredi} was a certain quantitative aspect that will not be discussed here, but the proof given there will turn out to be quite relevant for the present paper.

Generalizing Bourgain's result in another direction, Lyall and Magyar \cite{LM19:distgraphs} studied density theorems for the so-called \emph{distance graphs}. Informally speaking, these are graphs embedded in a Euclidean space that carry information about lengths of their edges. Certain non-degeneracy conditions are then needed in order to have meaningful results. We will neither give a precise definition of those concepts here, nor formulate the most general known result on this topic, which is \cite[Theorem~2]{LM19:distgraphs}. Instead, we will state its particular case when the graph is a tree, as this is the one that we are about to generalize later.

Take a tree $\mathcal{T}=(V,E)$ on a finite set of vertices $V$, having $E$ as the set of its edges. Suppose that we are also given a function $\ell\colon E\to(0,\infty)$, so that $\ell(e)$ is interpreted as the ``length'' of an edge $e\in E$. One could say that $\mathcal{T}$ equipped with $\ell$ is a \emph{distance tree}.
A special case of \cite[Theorem~2]{LM19:distgraphs} by Lyall and Magyar reads as follows:
\begin{itemize}
\item[] \emph{For every measurable set $A\subseteq\mathbb{R}^{2}$ satisfying $\overline{\delta}(A)>0$ there is a positive number $\lambda_0=\lambda_0(A,\mathcal{T},\ell)$ such that for each $\lambda\in[\lambda_0,\infty)$ the set $A$ contains a set of points $\{x_v:v\in V\}$ satisfying $|x_u-x_v|=\lambda\ell(e)$ whenever $e\in E$ is an edge connecting vertices $u,v\in V$.}
\end{itemize}
Interestingly, distance trees and their even more special cases, \emph{distance chains}, also play a prominent role in somewhat related problems \cite{BIT16:chains,IT19:trees}.
Distance trees are not rigid point configurations. One can still talk about their ``isometric copies'' within $A$ (and remain in line with the previous formulations), if one defines the concept of isometric distance graphs in an obvious way; see \cite{LM19:distgraphs}.

In formulations of all of these results the emphasis needs to be put on the fact that \emph{all} sufficiently large dilates of the configuration can be identified in the set $A$. Note that finding just \emph{any} dilate is trivial, since the Lebesgue density arguments identify sufficiently small dilates of any given finite configuration inside a set of positive measure. On the other hand, the existence of merely \emph{some} sufficiently large dilates can be deduced easily from Szemer\'{e}di's theorem \cite{S75:Szem} or its multidimensional version by Furstenberg and Katznelson \cite{FK78:msz}.

Let us also mention that there are many related papers that study lower-dimensional subsets $A$ of the Euclidean space \cite{BIT16:chains,CLP16:sparse,GIP17:necklaces,IL19:triangles,IT19:trees,GILP16:simplices,GIM18:thin,Y19:splx,IM20:simplices}, subsets $A$ of the multidimensional integer lattice \cite{M09:Zn,B18:Zn,LM15:Zn}, or patterns with arithmetic structure \cite{LP09:ap,Ke08:ap,HLP16:poly,S17:ap,FP18:large,CMP15:roth,DKR17:corner,DK18:products,DPZ19:rough,FGP19:roth,DK20:Szemeredi}. We do not discuss these types of results here.

\subsection{Statements of new results}
\label{subsec:newresults}
It is natural to start asking questions on generalizations of the above results to configurations that are not dilated equally in all directions. For instance, one might ask if a positive upper Banach density set necessarily contains copies of rectangles with sides $\lambda$ and $\lambda^2$ for all sufficiently large numbers $\lambda$. In the present paper we attempt to answer a few questions of this type, without insisting on formulating the most general possible results; see the comments in Subsection~\ref{subsec:comments} below. The presented proofs will reveal which of the techniques found in the literature allow modifications applicable to the anisotropic setting, and which objects from harmonic analysis appear along the way.

First, we turn to $n$-dimensional simplices in $\mathbb{R}^{n+1}$. Let us try to come up with a natural polynomial formulation of Bourgain's result. For instance, it does not make sense to look for triangles with sides $\lambda$, $\lambda^2$, and $\lambda^3$, because these three numbers fail to satisfy the triangle inequality for large $\lambda$. However, one can still be requiring that adjacent sides of the desired triangle are $\lambda$ and $\lambda^2$ units long and that the angle between them is fixed. Similarly, it can also be interesting to study right-angled simplices with perpendicular edges of lengths $\lambda, \lambda^2, \ldots, \lambda^n$, i.e., they are scaled according to the so-called moment curve $\lambda\mapsto(\lambda,\lambda^2,\ldots,\lambda^n)$.

In fact, polynomials will not play any role here and the approach will generalize naturally to power-type dilations with ratios of the form $\lambda^a b$, where $a,b>0$ are fixed parameters. In particular, $a$ does not need to be an integer. In other words, we will be working with \emph{anisotropic power-type dilations}
\begin{equation}\label{eq:anidilsm}
(x_1,\ldots,x_n) \mapsto (\lambda^{a_1} b_1 x_1,\ldots \lambda^{a_n} b_n x_n),
\end{equation}
but not necessarily in the standard coordinate system.
Here and in what follows, $n$ is a fixed positive integer, while
\begin{equation}\label{eq:asandbs}
a_1,a_2,\ldots,a_n, b_1,b_2,\ldots,b_n>0
\end{equation}
are fixed parameters.

Our first result is an anisotropic version of the above theorem of Bourgain \cite{B86:roth}.
Suppose that we are also given linearly independent unit vectors
\begin{equation}\label{eq:justus}
u_1, u_2, \ldots, u_n \in \mathbb{R}^{n},
\end{equation}
in addition to the numbers \eqref{eq:asandbs}. The idea is that these vectors determine directions of the edges (meeting at a single point) of the simplex with vertices $\Delta = \{\mathbf{0}, b_1 u_1, b_2 u_2, \ldots,$ $b_n u_n\}$.

\begin{theorem}\label{thm:simplices}
For every measurable set $A\subseteq\mathbb{R}^{n+1}$ satisfying $\overline{\delta}(A)>0$ there is a positive number $\lambda_0=\lambda_0(A,a_1,\ldots,a_n,b_1,\ldots,b_n,u_1,\ldots,u_n)$ such that for each $\lambda\in[\lambda_0,\infty)$ one can find a point $x\in\mathbb{R}^{n+1}$ and vectors $y_1,y_2,\ldots,y_n\in\mathbb{R}^{n+1}$ satisfying
\[ \{x,x+y_1,x+y_2,\ldots,x+y_n\} \subseteq A \]
and
\[ y_k\cdot y_l = \lambda^{a_k} b_k u_k\cdot \lambda^{a_l} b_l u_l \quad \text{for } k,l=1,2,\ldots,n. \]
In other words, for each $\lambda\in[\lambda_0,\infty)$ the set $A$ contains an isometric copy of
\[ \big\{\mathbf{0}, \lambda^{a_1} b_1 u_1, \lambda^{a_2} b_2 u_2, \ldots, \lambda^{a_n} b_n u_n\big\}. \]
\end{theorem}

A notable particular case of Theorem~\ref{thm:simplices} is obtained when the unit vectors \eqref{eq:justus} are mutually orthogonal, i.e., the simplex in question is right-angled. In this case the theorem simply guarantees the existence of mutually orthogonal vectors $y_k$ with lengths $\lambda^{a_k} b_k$ such that a translate of $\{\mathbf{0},y_1,\ldots,y_n\}$ is contained in $A$.

Now we turn to $n$-dimensional rectangular boxes in $\mathbb{R}^{2n}$.
Our second result is an anisotropic generalization of the above theorem of Lyall and Magyar \cite{LM19:hypergraphs}.

\begin{theorem}\label{thm:boxes}
For every measurable set $A\subseteq(\mathbb{R}^{2})^n$ satisfying $\overline{\delta}(A)>0$ there is a positive number $\lambda_0=\lambda_0(A,a_1,\ldots,a_n,b_1,\ldots,b_n)$ such that for each $\lambda\in[\lambda_0,\infty)$ one can find
$x_1,\ldots,x_n, y_1,\ldots,y_n \in\mathbb{R}^2$ satisfying
\[ \big\{ (x_1 + r_1 y_1, x_2 + r_2 y_2, \ldots, x_n + r_n y_n) : (r_1,\ldots,r_n)\in\{0,1\}^n \big\} \subseteq A \]
and
\[ |y_k| = \lambda^{a_k} b_k\quad \text{for } k=1,2,\ldots,n. \]
In other words, for each $\lambda\in[\lambda_0,\infty)$ the set $A$ contains an isometric copy of
\[ \{0,\lambda^{a_1} b_1\} \times \{0,\lambda^{a_2} b_2\} \times \cdots\times \{0,\lambda^{a_n} b_n\} \subset \mathbb{R}^n \]
with sides parallel to the distinguished $2$-dimensional coordinate planes.
\end{theorem}

At the first sight it might appear that Theorem~\ref{thm:boxes} generalizes the ``right-angled case'' of Theorem~\ref{thm:simplices}, since vertices of a right-angled simplex clearly form a subset of the set of vertices of an appropriate rectangular box. A subtle distinction is that Theorem~\ref{thm:simplices} already holds in the $(n+1)$-dimensional Euclidean space, but we allowed the simplex to rotate in all possible directions. As opposed to that, Theorem~\ref{thm:boxes} is easily seen to fail in less than $2n$ dimensions, because in it we consider only rotations coming from $n$ coordinate planes of the splitting $(\mathbb{R}^{2})^n=\mathbb{R}^2\times\mathbb{R}^2\times\cdots\times\mathbb{R}^2$.
A $(2n-1)$-dimensional counterexample is the set $A\subseteq \mathbb{R}\times\mathbb{R}^2\times\cdots\times\mathbb{R}^2$ obtained by restricting the first coordinate to $\bigcup_{m\in\mathbb{Z}}[m-1/10,m+1/10]$.
On the other hand, the corresponding result for arbitrarily rotated rectangular boxes in less than $2n$ dimensions has not been either proved or disproved at the time of writing, even in the case of a cube, i.e., when all parameters from \eqref{eq:asandbs} are equal to $1$.

Finally, we give an anisotropic generalization of the aforementioned result of Lyall and Magyar on distance trees \cite{LM19:distgraphs}.
Let $\mathcal{T}=(V,E)$ be a finite tree with vertices $V$ and edges $E$. It is convenient to identify the set of edges $E$ with $\{1,2,\ldots,n\}$ and the number of vertices is then equal to $n+1$. This way the parameters $a_k,b_k$ from \eqref{eq:asandbs} are associated to the edges $k\in E$ of the tree. We no longer need to mention any length function $\ell$, as the assignment $k\mapsto a_k,b_k$ gives rise to an even more complex structure. However, if we want to have some length function defined explicitly, then we can set $\ell(k):=b_k$ for each edge $k$.

\begin{theorem}\label{thm:trees}
For every measurable set $A\subseteq\mathbb{R}^{2}$ satisfying $\overline{\delta}(A)>0$ there is a positive number $\lambda_0=\lambda_0(A,\mathcal{T},a_1,\ldots,a_n,b_1,\ldots,b_n)$ such that for each $\lambda\in[\lambda_0,\infty)$ one can find a set of points
\[ \{ x_v : v\in V \} \subseteq A \]
satisfying
\[ |x_u-x_v| = \lambda^{a_k} b_k\quad \text{for each edge } k\in E \text{ joining vertices } u,v\in V. \]
In other words, for each $\lambda\in[\lambda_0,\infty)$ the set $A$ contains an embedding of the distance tree combinatorially isomorphic to $\mathcal{T}$ and having the numbers $\ell(k)=\lambda^{a_k} b_k$ as lengths of its edges.
\end{theorem}

Note that Theorem~\ref{thm:trees} is placed in two dimensions only. If we disregarded its dimensional sharpness, the particular case $a_1=\cdots=a_n$, $b_1=\cdots=b_n$ of Theorem~\ref{thm:trees} would be a consequence of Theorem~\ref{thm:boxes}, because each tree is easily seen to be a subgraph of the hypercube graph in a sufficiently large dimension.

Proofs of Theorems~\ref{thm:simplices}--\ref{thm:trees} will rely on a few relatively known ideas from the harmonic analysis, modulo a general approach described in Subsection~\ref{subsec:approach}.
Therefore, the main contribution of the present paper lies simply in recollecting, selecting, and reapplying those ideas to the above combinatorial problems.
Further connections between the harmonic analysis and the combinatorics of the Euclidean space will be discussed at the end of the paper, in Subsection~\ref{subsec:anisosingint}.

\subsection{General scheme of the approach}
\label{subsec:approach}
The proofs of Theorems~\ref{thm:simplices}, \ref{thm:boxes}, and \ref{thm:trees} will be presented in Sections~\ref{sec:simplices}, \ref{sec:boxes}, and \ref{sec:trees}, respectively. Here we only discuss the general outline.

We will fit the proofs to the scheme that we are about to describe. The approach is an abstraction of the method stemming from the work of Bourgain \cite{B86:roth} and first used by Cook, Magyar, and Pramanik \cite{CMP15:roth} in a way that emphasizes the role of estimates for multilinear singular integrals or similar objects. Its variant was named the \emph{largeness--smoothness multiscale approach} by Durcik and the present author \cite{DK20:Szemeredi}. Limitations of the method are essentially only the limitations within the field of multilinear harmonic analysis, which has seen vast and rapid development over the last few decades. The approach has already been reused several times after \cite{CMP15:roth} (see the papers \cite{DKR17:corner,DK18:products,DK20:Szemeredi}, which study arithmetic progressions and related configurations), but here we want to point out that the method is also effective for many geometric configurations without any algebraic structure.

For each of the studied problems we will define a \emph{counting form} $\mathcal{N}^{0}_{\lambda}$ that identifies the configuration associated with the parameter $\lambda>0$. In order to prove the claim it is sufficient to show that $\mathcal{N}^{0}_{\lambda}$ is positive for all $\lambda$ that are sufficiently large depending on the set $A$. Besides $\lambda$, which can be thought of as a \emph{scale of largeness}, there will be another scale $0<\varepsilon\leq 1$, interpreted as a \emph{scale of smoothness}. A two-parameter family of counting forms $\mathcal{N}^{\varepsilon}_{\lambda}$ will recover $\mathcal{N}^{0}_{\lambda}$ in the limit as $\varepsilon\to0$. The reason for ``smoothing'' or ``blurring out'' up to scale $\varepsilon$ is that the smoother configuration can be identified with a more direct counting argument. Thus, the method starts by decomposing $\mathcal{N}^{0}_{\lambda}$ as
\begin{equation}\label{eq:basicsplitting}
\mathcal{N}^{1}_{\lambda} + \big(\mathcal{N}^{\varepsilon}_{\lambda}-\mathcal{N}^{1}_{\lambda}\big) + \big(\mathcal{N}^{0}_{\lambda}-\mathcal{N}^{\varepsilon}_{\lambda}\big).
\end{equation}

The smoothest term $\mathcal{N}^{1}_{\lambda}$ can be thought of as the \emph{structured part} and its lower bound is a simpler problem.
Typically one needs to zoom the picture to scale $\lambda$ and then perform a direct counting argument. In this paper we need to handle mismatched scales $\lambda^{a_1},\ldots,\lambda^{a_n}$ simultaneously, which is a novel complication in relation with simplices, boxes, or trees, but it has already appeared in the work of Bourgain on polynomial three-term progressions \cite{B88:nonlinroth}.
The term $\mathcal{N}^{1}_{\lambda}$ is a reason why we will first show several estimates for general filtrations of a fixed probability space in Subsection~\ref{subsec:filtrations}.

The third term $\mathcal{N}^{0}_{\lambda}-\mathcal{N}^{\varepsilon}_{\lambda}$ is interpreted as the \emph{uniform part} and some oscillatory phenomenon should guarantee that it converges to $0$ uniformly in $\lambda$ as $\varepsilon\to0$. Quite often and also in this paper, the only oscillatory estimate needed is the decay of the Fourier transform of a spherical measure; see \eqref{eq:circledecay} and \eqref{eq:subcircledecay} below. Thus, one can fix a sufficiently small $\varepsilon>0$ such that the uniform part is always dominated by the structured part.

The middle term $\mathcal{N}^{\varepsilon}_{\lambda}-\mathcal{N}^{1}_{\lambda}$ is the \emph{error part}. It cannot be efficiently estimated for a fixed value of $\lambda$, so we rather attempt to control it ``on the average'' for sufficiently many scales $\lambda_1<\lambda_2<\cdots<\lambda_J$ satisfying, say, $\lambda_{j+1}\geq 2\lambda_j$ for each $j$. More precisely, sums of the form
\begin{equation}\label{eq:mltiscaleqnt}
\sum_{j=1}^{J} \big|\mathcal{N}^{\varepsilon}_{\lambda_j}-\mathcal{N}^{1}_{\lambda_j}\big|
\end{equation}
for lacunary scales $\lambda_j$ are shown to satisfy a bound that is allowed to blow up as $\varepsilon\to0$, but it is at the same nontrivial in the total number of scales $J$, i.e., it grows like $o(J)$ as $J\to\infty$. When $J$ is sufficiently large, pigeonholing guarantees that at least one of the individual errors $|\mathcal{N}^{\varepsilon}_{\lambda_j}-\mathcal{N}^{1}_{\lambda_j}|$ is sufficiently small.
It is precisely the multiscale quantity \eqref{eq:mltiscaleqnt} that resembles a certain multilinear singular integral form.
Bounds for \eqref{eq:mltiscaleqnt} are shown using certain ``cancellation'' between different scales $\lambda_j$. This can be done using techniques from multilinear harmonic analysis, by treating \eqref{eq:mltiscaleqnt} as a multisublinear integral operator.
This is a route that we follow in the present paper, except that we clean up the proofs by reducing the aforementioned operator bounds merely to several Gaussian identities and estimates that we first establish in Subsection~\ref{subsec:Gaussianiden}. To a large extent these have already appeared in \cite{DK20:Szemeredi}.

Finally, the actual result is shown by contradiction: assuming that the set $A$ contains no copies of the desired configuration associated with a lacunary sequence of parameters $\lambda_1<\lambda_2<\cdots$.
By choosing $\varepsilon>0$ sufficiently small and then choosing $J$ sufficiently large, we can guarantee that $|\mathcal{N}^{\varepsilon}_{\lambda_j}-\mathcal{N}^{1}_{\lambda_j}|$ and $|\mathcal{N}^{0}_{\lambda_j}-\mathcal{N}^{\varepsilon}_{\lambda_j}|$ are both dominated by $\mathcal{N}^{1}_{\lambda_j}$ for at least one index $j$. This implies $\mathcal{N}^{0}_{\lambda_j}>0$, which contradicts our hypothesis that $A$ does not contain the desired configuration of size $\lambda_j$.
A small technicality is that it is more convenient to work with a localized version $B$ of the given set $A$.
Details of the method applied to specific configurations can be found, for instance, in \cite[Section~2]{CMP15:roth} or \cite[Section~3]{DKR17:corner}. We will also be completely rigorous about these details in Sections~\ref{sec:simplices}--\ref{sec:trees}.

Already the pioneering work of Bourgain \cite{B86:roth} used a particular case of the above scheme of proof. The main novelty introduced by Cook, Magyar, and Pramanik \cite{CMP15:roth} is that the smoother version $\mathcal{N}^{\varepsilon}_{\lambda}$ of the counting form $\mathcal{N}^{0}_{\lambda}$ need not be obtained by smoothing the input functions (see Sections~\ref{sec:boxes} and \ref{sec:trees}), even though this is one legitimate possibility (see Section~\ref{sec:simplices}). Admittedly, this general scheme was primarily devised for studying ``more singular configurations,'' such as arithmetic progressions; see \cite{CMP15:roth,DK20:Szemeredi,DK18:products,DKR17:corner}. It can be an overkill in relation with simplices or boxes, as the papers \cite{B86:roth,LM19:distgraphs,LM16:prod} proceed by following a different philosophy. However, the present work does benefit from the power and flexibility of the general largeness--smoothness multiscale approach. Namely, we will define $\mathcal{N}^{\varepsilon}_{\lambda}$ using the heat flow --- the motivation comes from \cite[Section~7]{DK20:Szemeredi}, while \cite[Sections~3--6]{DK20:Szemeredi} also employ a similar time-space dynamics. That way, we will make use of a simple fact that the heat equation remains essentially the same after a power-type change of the time variable; compare Formulae \eqref{eq:heateq} and \eqref{eq:heateqgen} below.

As we said, Gaussians play a prominent role throughout the paper, so Subsection~\ref{subsec:Gaussianiden} will recall their properties needed in the proofs. Most notable ones will be the heat equation \eqref{eq:heateq}, estimates \eqref{eq:GaussianbySchwartz} and \eqref{eq:SchwartzbyGaussian}, and identities \eqref{eq:idt1} and \eqref{eq:infortheta} below.
They allow us to easily estimate convolutions with general probability measures both from below (which will be needed in Subsections~\ref{subsec:simplexstructured}, \ref{subsec:boxesstructured}, \ref{subsec:structrees}) and from above (which will be needed in Subsections~\ref{subsec:simplexerror}, \ref{subsec:boxeserror}, \ref{subsec:errtrees}).
It is quite likely that other semigroup structures work as well, at least in Section~\ref{sec:simplices} as Bourgain \cite{B86:roth} used the Poisson kernel for simplices. In any case, the present paper tries to advertise Gaussians as convenient mollifiers for the problems studied here.

\subsection{Organization of the paper}
Let us explain shortly how the rest of the paper is organized.
Section~\ref{sec:notation} discusses the notation used throughout the paper. It also recalls a few basic notions from the Fourier analysis, proves a couple of identities concerning Gaussian functions, and shows a few inequalities for conditional expectations on a general probability space.
Section~\ref{sec:simplices} establishes Theorem~\ref{thm:simplices} on anisotropic simplices, Section~\ref{sec:boxes} establishes Theorem~\ref{thm:boxes} on anisotropic rectangular boxes, while Section~\ref{sec:trees} proves Theorem~\ref{thm:trees} on anisotropic trees.
Each of these sections is divided further into three subsections that respectively handle structured, error, and uniform terms from the basic splitting \eqref{eq:basicsplitting}.
Finally, Section~\ref{sec:closing} discusses anisotropic multilinear singular integral operators that are motivated by the above combinatorial problems.
It also comments on possible generalizations of the results and limitations of the approach.

\section{Notation and preliminaries}
\label{sec:notation}

\subsection{Basic notation}
Let $A$ and $B$ be two nonnegative quantities. We write $A\lesssim B$ and $B\gtrsim A$ if $A\leq C B$ holds for some (unimportant) finite positive constant $C$. We write $A\sim B$ if $c B\leq A\leq C B$ holds for finite positive constants $c$ and $C$. Throughout the paper it will be understood that any of these constants $c,C$ are allowed to depend on the dimension of the ambient Euclidean space, the number $n$, the parameters from \eqref{eq:asandbs}, and (in Section~\ref{sec:simplices}) also on the unit vectors from \eqref{eq:justus}, but are independent of all other parameters or variables.

The open Euclidean ball with radius $r$ centered at $x$ will be denoted $\textup{B}(x,r)$.
The (standard) \emph{inner product} and the \emph{Euclidean norm} on $\mathbb{R}^d$ will be written as $(x,y)\mapsto x\cdot y$ and $x\mapsto |x|$, respectively.
The distance from a point $x\in\mathbb{R}^d$ to a set $S\subseteq\mathbb{R}^d$ will be denoted $\mathop{\textup{dist}}(x,S)$.
The linear span of a set of vectors $S\subseteq \mathbb{R}^d$ will be written as $\mathop{\textup{span}}(S)$.
We will write $\mathbbm{1}_A$ for the \emph{indicator function} of a set $A\subseteq\mathbb{R}^d$.
The \emph{floor} function is denoted $x\mapsto\lfloor x\rfloor$, i.e., $\lfloor x\rfloor$ is the largest integer less than or equal to $x\in\mathbb{R}$.
The complex \emph{imaginary unit} will be written as  $\mathbbm{i}$.
The logarithm function will be written ``$\log$'' and it will be understood that its base is the number $e$.

We will always specify the measure with respect to which the integrals are evaluated, unless we are working with the Lebesgue measure. Similarly, we will simply write $|A|$ for the Lebesgue measure of $A$.
By $\fint_{A}f$ we denote the \emph{average} of a locally integrable complex function $f$ over a bounded measurable set $A\subseteq\mathbb{R}^d$.
We will write $f\mapsto\|f\|_{\textup{L}^p}$ for the norm of $\textup{L}^p(\mathbb{R}^d)$, $p\in[1,\infty]$, and $(f,g)\mapsto\langle f,g\rangle_{\textup{L}^2}$ for the inner product in $\textup{L}^2(\mathbb{R}^d)$.

Now we come to \emph{dilates} and \emph{convolutions} of functions and measures.
For an integrable function $f\colon\mathbb{R}^d\to\mathbb{C}$ and a number $\lambda\in\mathbb{R}\setminus\{0\}$ we define
\[ f_{\lambda}(x) := |\lambda|^{-d} f(\lambda^{-1}x) \quad \text{for }x\in\mathbb{R}^d. \]
More generally, for a finite measure $\nu$ on Borel subsets $A\subseteq\mathbb{R}^d$ we write
\[ \nu_{\lambda}(A) := \nu(\lambda^{-1}A)
= \nu(\{\lambda^{-1}x : x\in A\}). \]
Note that the normalizations of $f_\lambda$ and $\nu_\lambda$ are consistent with each other: if $\nu$ happens to be an absolutely continuous measure with density $f$, then $\nu_\lambda$ will have $f_\lambda$ for its density.
In fact, for a bounded measurable function $h\colon\mathbb{R}^d\to\mathbb{C}$ we have
\begin{align*}
\int_{\mathbb{R}^d} h(x) \,\textup{d}\nu_{\lambda}(x) & = \int_{\mathbb{R}^d} h(\lambda x) \,\textup{d}\nu(x), \\
\int_{\mathbb{R}^d} h(x) f_{\lambda}(x) \,\textup{d}x & = \int_{\mathbb{R}^d} h(\lambda x) f(x) \,\textup{d}x.
\end{align*}
If $g\colon\mathbb{R}^d\to\mathbb{C}$ is another integrable function, then it makes sense to define
\[ (f\ast g)(x) := \int_{\mathbb{R}^d} f(y) g(x-y) \,\textup{d}y \quad \text{for a.e.\@ }x\in\mathbb{R}^d. \]
It is well known that the operation of convolution is commutative.
More generally, the convolution of a finite measure $\nu$ and an integrable function $g$ is defined as
\[ (\nu\ast g)(x) := \int_{\mathbb{R}^d} g(x-y) \,\textup{d}\nu(y) \quad \text{for a.e.\@ }x\in\mathbb{R}^d. \]
Even at this level of generality the operation of convolution is associative, i.e., $(\nu\ast f)\ast g=\nu\ast (f\ast g)$, while the Dirac delta measure at the origin, denoted $\delta_{\mathbf{0}}$, serves as the identity element.

The \emph{Fourier transform} of an integrable function $f\colon\mathbb{R}^d\to\mathbb{C}$ is $\widehat{f}\colon\mathbb{R}^d\to\mathbb{C}$ defined as
\[ \widehat{f}(\xi) := \int_{\mathbb{R}^d} f(x) e^{-2\pi \mathbbm{i} x\cdot\xi} \,\textup{d}x \quad \text{for }\xi\in\mathbb{R}^d, \]
while the Fourier transform of a finite Borel measure $\nu$ is $\widehat{\nu}\colon\mathbb{R}^d\to\mathbb{C}$ defined as
\[ \widehat{\nu}(\xi) := \int_{\mathbb{R}^d} e^{-2\pi \mathbbm{i} x\cdot\xi} \,\textup{d}\nu(x) \quad \text{for }\xi\in\mathbb{R}^d. \]
Basic properties of the Fourier transform can be found in any introductory textbook on the harmonic analysis, such as \cite{SW71:book}.
For instance, it is useful to know that, for any $\lambda\in\mathbb{R}\setminus\{0\}$ and $f,\nu$ as above, one has
\[ \widehat{f_\lambda}(\xi) = \widehat{f}(\lambda \xi), \quad
\widehat{\nu_\lambda}(\xi) = \widehat{\nu}(\lambda \xi) \quad \text{for } \xi\in\mathbb{R}^d. \]

Important instances of Borel measures are \emph{spherical measures}. If $\sigma$ is the normalized surface measure of the $(d-1)$-dimensional standard unit sphere $\mathbb{S}^{d-1}$ in $\mathbb{R}^{d}$, $d\geq 2$, then its Fourier transform satisfies the well-known decay
\begin{equation}\label{eq:circledecay}
|\widehat{\sigma}(\xi)| \lesssim \min\big\{1,|\xi|^{-(d-1)/2}\big\} \quad \text{for } \xi\in\mathbb{R}^d;
\end{equation}
see \cite[Subsection~VIII.5.B]{St93:book}.
If $\sigma$ is a normalized $(k-1)$-dimensional spherical measure supported on a sphere of radius $r>0$ in a $k$-dimensional plane in $\mathbb{R}^{d}$ orthogonal to some $(d-k)$-dimensional linear subspace $H\subset\mathbb{R}^{d}$, then
\begin{equation}\label{eq:subcircledecay}
|\widehat{\sigma}(\xi)| \lesssim \min\big\{1,\big(r\mathop{\textup{dist}}(\xi,H)\big)^{-(k-1)/2}\big\} \quad \text{for } \xi\in\mathbb{R}^d.
\end{equation}

\subsection{Gaussian identities}
\label{subsec:Gaussianiden}
We write $\mathbbm{g}$ for the \emph{standard $d$-dimensional Gaussian function},
\[ \mathbbm{g}\colon\mathbb{R}^d\to[0,\infty), \quad \mathbbm{g}(x) := e^{-\pi|x|^2}. \]
We also reserve special letters for its partial derivatives,
\[ \mathbbm{h}^{(l)} := \partial_l \mathbbm{g} \quad  \text{for } l=1,2,\ldots,d \]
and for its Laplacian
\[ \mathbbm{k} := \Delta\mathbbm{g}. \]
Basic properties of the Fourier transform easily yield
\[ \widehat{\mathbbm{g}}(\xi) = \mathbbm{g}(\xi), \]
\[ \widehat{\mathbbm{h}^{(l)}}(\xi) = 2\pi\mathbbm{i}\xi_l \,\mathbbm{g}(\xi) \quad  \text{for } l=1,2,\ldots,d, \]
and
\[ \widehat{\mathbbm{k}}(\xi) = - 4 \pi^2 |\xi|^2 e^{-\pi |\xi|^2}, \]
where $\xi=(\xi_1,\xi_2,\ldots,\xi_d)\in\mathbb{R}^d$ is arbitrary. These formulae (and the fact that the Fourier transform interchanges convolutions and pointwise products) imply the following convolution identities:
\begin{equation}\label{eq:cvid1}
\mathbbm{g}_{\alpha} \ast \mathbbm{g}_{\beta} = \mathbbm{g}_{\sqrt{\alpha^2+\beta^2}},
\end{equation}
\begin{equation}\label{eq:cvid2}
\sum_{l=1}^{d} \mathbbm{h}^{(l)}_{\alpha} \ast \mathbbm{h}^{(l)}_{\beta} = \frac{\alpha\beta}{\alpha^2+\beta^2} \mathbbm{k}_{\sqrt{\alpha^2+\beta^2}},
\end{equation}
and
\begin{equation}\label{eq:cvid3}
\mathbbm{k}_{\alpha} \ast \mathbbm{g}_{\beta} = \frac{\alpha^2}{\alpha^2+\beta^2} \mathbbm{k}_{\sqrt{\alpha^2+\beta^2}}
\end{equation}
for any $\alpha,\beta\in(0,\infty)$.
Identities \eqref{eq:cvid1}--\eqref{eq:cvid3} will be useful for splitting $\mathbbm{g}_t$ or $\mathbbm{k}_t$ into a convolution of one term with a desired scale and a uniquely determined remaining term; see Subsections~\ref{subsec:simplexerror} and \ref{subsec:boxeserror}.

The above Gaussian functions are easily seen to satisfy the \emph{heat equation}:
\begin{equation}\label{eq:heateq}
\frac{\partial}{\partial t} \big( \mathbbm{g}_{t}(x) \big) = \frac{1}{2\pi t} \mathbbm{k}_{t}(x)
\end{equation}
on $(t,x)\in (0,\infty)\times\mathbb{R}^d$.
By a simple chain rule, \eqref{eq:heateq} generalizes to
\begin{equation}\label{eq:heateqgen}
\frac{\partial}{\partial t} \big( \mathbbm{g}_{t^a b}(x) \big) = \frac{a}{2\pi t} \mathbbm{k}_{t^a b}(x),
\end{equation}
where $a,b\in(0,\infty)$ can be arbitrary.
The heat equation will govern the smoothing dynamics; see the beginning of Subsection~\ref{subsec:simplexerror} and the beginning of Subsection~\ref{subsec:boxeserror}.

On the one hand, since Gaussian tails decay faster than any polynomial, we trivially have
\begin{equation}\label{eq:GaussianbySchwartz}
\mathbbm{g}(x)\lesssim (1+|x|)^{-d-1}, \quad |\mathbbm{h}^{(l)}(x)|\lesssim (1+|x|)^{-d-1}
\end{equation}
for $x\in\mathbb{R}^d$.
On the other hand,
\begin{equation}\label{eq:SchwartzbyGaussian}
(1+|x|)^{-d-1} \lesssim \int_1^\infty \mathbbm{g}_{\gamma}(x) \,\frac{\textup{d}\gamma}{\gamma^{2}}
\end{equation}
for $x\in\mathbb{R}^d$. In words, Schwartz tails are dominated by a superposition of dilated Gaussians. Formula \eqref{eq:SchwartzbyGaussian} was first used in a similar context by Durcik \cite{D15:L4}. It can be shown easily by investigating the asymptotic behavior of the right hand side as $|x|\to\infty$; see the details in \cite[Section~3]{D15:L4} or \cite[Section~3]{DKST19:NVEA}.
A combination of \eqref{eq:GaussianbySchwartz} and \eqref{eq:SchwartzbyGaussian} will be used to bound a convolution of a Gaussian and an arbitrary probability measure with bounded support by a superposition of (non-translated) Gaussians; see the computations leading to \eqref{eq:Gaussdomh} and \eqref{eq:Gaussdomg} below.

We claim a simple identity:
\begin{equation}\label{eq:idt1}
\int_{0}^{\infty} \sum_{l=1}^{d} \big\|f\ast\mathbbm{h}^{(l)}_{s^{a}b}\big\|_{\textup{L}^2}^2 \,\frac{\textup{d}s}{s} = \frac{\pi}{a} \|f\|_{\textup{L}^2}^2
\end{equation}
for any compactly supported $f\in\textup{L}^2(\mathbb{R}^d)$ and $a,b\in(0,\infty)$.
Indeed, using \eqref{eq:cvid2}, \eqref{eq:heateqgen}, and \eqref{eq:cvid1}, respectively, the left hand side of \eqref{eq:idt1} can be rewritten as
\begin{align*}
& -\lim_{\substack{\alpha\to0^+\\ \beta\to\infty}} \int_{\alpha}^{\beta} \sum_{l=1}^{d} \big\langle f\ast\mathbbm{h}^{(l)}_{s^{a}b}\ast\mathbbm{h}^{(l)}_{s^{a}b}, f \big\rangle_{\textup{L}^2} \,\frac{\textup{d}s}{s} \\
& = \frac{\pi}{a} \lim_{\substack{\alpha\to0^+\\ \beta\to\infty}} \Big( \big\langle f\ast(\mathbbm{g}_{2^{1/2}\alpha^{a}b}-\mathbbm{g}_{2^{1/2}\beta^{a}b}), f \big\rangle_{\textup{L}^2} \Big) \\
& = \frac{\pi}{a} \Big( \lim_{\alpha\to0^+} \big\|f\ast\mathbbm{g}_{\alpha^{a}b}\big\|_{\textup{L}^2}^2
- \lim_{\beta\to\infty} \big\|f\ast\mathbbm{g}_{\beta^{a}b}\big\|_{\textup{L}^2}^2 \Big)
= \frac{\pi}{a} \|f\|_{\textup{L}^2}^2.
\end{align*}

Next, for fixed parameters \eqref{eq:asandbs}, for any choice of $\gamma_1,\ldots,\gamma_n\in(0,\infty)$, and for a compactly supported real-valued $f\in\textup{L}^{2^n}((\mathbb{R}^d)^n)$ we define
\begin{align}
\Theta^{n,m}_{\gamma_1,\ldots,\gamma_n}(f) :=
- \int_{0}^{\infty} \int_{(\mathbb{R}^d)^{2n}} & \mathcal{F}(\mathbf{x})
\, \mathbbm{k}_{s^{a_m}b_m\gamma_m}(x_m^0-x_m^1) \nonumber \\ & \times \Big(\prod_{\substack{1\leq k\leq n\\ k\neq m}}\mathbbm{g}_{s^{a_k}b_k\gamma_k}(x_k^0-x_k^1)\Big)
\,\textup{d}\mathbf{x} \,\frac{\textup{d}s}{s}. \label{eq:defofthetaform}
\end{align}
Here we denote
\begin{equation}\label{eq:boxauxdef1}
\mathcal{F}(\mathbf{x}) := \prod_{(r_1,\ldots,r_n)\in\{0,1\}^n} f(x_1^{r_1}, \ldots, x_n^{r_n}),
\end{equation}
so that this is a function of $\mathbf{x}=(x_1^0,x_1^1,\ldots,x_n^0,x_n^1)\in(\mathbb{R}^d)^{2n}$, and write formally
\begin{equation}
\textup{d}\mathbf{x} := \textup{d}x_1^0\,\textup{d}x_1^1 \,\textup{d}x_2^0\,\textup{d}x_2^1 \cdots \textup{d}x_n^0\,\textup{d}x_n^1.
\end{equation}
We will also denote
\begin{equation}\label{eq:boxauxdef2}
\mathcal{F}^{(m)}(x) := \prod_{(r_1,\ldots,r_{m-1},r_{m+1}\ldots,r_n)\in\{0,1\}^{n-1}}
f(x_1^{r_1}, \ldots, x_{m-1}^{r_{m-1}}, x, x_{m+1}^{r_{m+1}}, \ldots, x_{n}^{r_{n}})
\end{equation}
for $x\in\mathbb{R}^d$ and $m\in\{1,\ldots,n\}$, keeping in mind that $\mathcal{F}^{(m)}(x)$ also depends on other variables than just $x$.

Formula \eqref{eq:cvid2} allows us to rewrite
\begin{align*}
\Theta^{n,m}_{\gamma_1,\ldots,\gamma_n}(f)
= 2 \sum_{l=1}^{d} \int_{0}^{\infty} \int_{(\mathbb{R}^d)^{2n-2}} & \big\| \mathcal{F}^{(m)}\ast \mathbbm{h}_{2^{-1/2}s^{a_m}b_m\gamma_m}^{(l)} \big\|_{\textup{L}^2(\mathbb{R}^d)}^2 \\
& \times \Big(\prod_{\substack{1\leq k\leq n\\ k\neq m}}\mathbbm{g}_{s^{a_k}b_k\gamma_k}(x_k^0-x_k^1)\,\textup{d}x_k^0\,\textup{d}x_k^1\Big) \,\frac{\textup{d}s}{s},
\end{align*}
which reveals that $\Theta^{n,m}_{\gamma_1,\ldots,\gamma_n}(f)$ is nonnegative and well-defined.
This time we claim the identity
\begin{equation}\label{eq:infortheta}
\sum_{m=1}^{n} a_m \Theta^{n,m}_{\gamma_1,\ldots,\gamma_n}(f) = 2\pi \|f\|_{\textup{L}^{2^n}(\mathbb{R}^d)}^{2^n}.
\end{equation}
First, by the product rule for differentiation and the generalized heat equation \eqref{eq:heateqgen} we can write
\begin{equation}\label{eq:theprodrule}
\frac{\partial}{\partial s} \Big(\prod_{k=1}^{n}\mathbbm{g}_{s^{a_k}b_k\gamma_k}(y_k)\Big)
= \sum_{m=1}^{n} \frac{a_m}{2\pi s} \mathbbm{k}_{s^{a_m}b_m\gamma_m}(y_m) \Big(\prod_{\substack{1\leq k\leq n\\ k\neq m}}\mathbbm{g}_{s^{a_k}b_k\gamma_k}(y_k)\Big).
\end{equation}
A consequence of the last display and the fundamental theorem of calculus is
\begin{align*}
\sum_{m=1}^{n} a_m \Theta^{n,m}_{\gamma_1,\ldots,\gamma_n}(f)
& = 2\pi \lim_{\alpha\to0^+} \int_{(\mathbb{R}^d)^{2n}} \mathcal{F}(\mathbf{x}) \Big(\prod_{k=1}^{n}\mathbbm{g}_{\alpha^{a_k}b_k\gamma_k}(x_k^0-x_k^1)\Big) \,\textup{d}\mathbf{x} \\
& \quad - 2\pi \lim_{\beta\to\infty} \int_{(\mathbb{R}^d)^{2n}} \mathcal{F}(\mathbf{x}) \Big(\prod_{k=1}^{n}\mathbbm{g}_{\beta^{a_k}b_k\gamma_k}(x_k^0-x_k^1)\Big) \,\textup{d}\mathbf{x}.
\end{align*}
The first limit above equals (allowing a slighly informal usage of $\delta_{\mathbf{0}}$):
\begin{align*}
& \int_{(\mathbb{R}^d)^{2n}} \mathcal{F}(\mathbf{x}) \Big(\prod_{k=1}^{n}\delta_{\mathbf{0}}(x_k^0-x_k^1)\Big) \,\textup{d}\mathbf{x} \\
& = \int_{(\mathbb{R}^d)^{n}} f(x_1,\ldots,x_n)^{2^n} \,\textup{d}x_1 \cdots \textup{d}x_n
= \|f\|_{\textup{L}^{2^n}(\mathbb{R}^d)}^{2^n},
\end{align*}
while the second one is $0$. This proves \eqref{eq:infortheta}.

Every summand on the left hand side of \eqref{eq:infortheta} is nonnegative, so the $\Theta^{n,m}_{\gamma_1,\ldots,\gamma_n}(f)$ is clearly bounded by $(2\pi/a_m)\|f\|_{\textup{L}^{2^n}(\mathbb{R}^d)}^{2^n}$.
Note that this bound is uniform in the parameters $\gamma_1,\ldots,\gamma_n$.
Identities \eqref{eq:idt1} and \eqref{eq:infortheta} will bound multiscale expressions coming from the study of the error parts in the decomposition \eqref{eq:basicsplitting}.
They can be viewed, respectively, as cheap substitutes for square function estimates and bounds for entangled multilinear singular integrals (mentioned in Subsection~\ref{subsec:anisosingint}).

\subsection{Conditional expectations}
\label{subsec:filtrations}
A few concepts from probability theory will be useful in the proofs below. Even though the dyadic setting would be sufficient, the general probabilistic notation makes arguments elegant and concise.

We are working in a fixed probability space $(\Omega,\mathcal{F},\mathbb{P})$. The associated \emph{expectation} is the operator $\mathbb{E}\colon f\mapsto\int_{\Omega} f\,\textup{d}\mathbb{P}$ defined on $\textup{L}^1(\Omega,\mathcal{F},\mathbb{P})$. Moreover, for any $\sigma$-algebra $\mathcal{G}\subseteq\mathcal{F}$ one can construct the operator of \emph{conditional expectation} with respect to $\mathcal{G}$ as the map
\[ \textup{L}^1(\Omega,\mathcal{F},\mathbb{P}) \to \textup{L}^1(\Omega,\mathcal{G},\mathbb{P}), \quad f\mapsto \mathbb{E}(f|\mathcal{G}) \]
such that
\[ \int_{A} \mathbb{E}(f|\mathcal{G})\,\textup{d}\mathbb{P} = \int_{A} f\,\textup{d}\mathbb{P} \]
for every $f\in \textup{L}^1(\Omega,\mathcal{F},\mathbb{P})$ and every $A\in\mathcal{G}$. The proof of its existence (and uniqueness) can be found in many textbooks on introductory probability theory, such as the one by Durrett \cite{Du19:book}. If $\mathcal{G}$ is generated by a finite partition of $\Omega$ into atoms and if, among these atoms, $A_1,\ldots,A_N$ have nonzero probability, then we have a simple formula
\begin{equation}\label{eq:ce1}
\mathbb{E}(f|\mathcal{G}) = \sum_{n=1}^{N} \Big(\frac{1}{\mathbb{P}(A_n)}\int_{A_n}f\,\textup{d}\mathbb{P}\Big) \mathbbm{1}_{A_n} \quad \textup{a.s.}
\end{equation}
In words, on each of the atoms $\mathbb{E}(f|\mathcal{G})$ is constantly equal to the average value of $f$ over that atom.
As a very special case of only one atom, we get
\begin{equation}\label{eq:ce2}
\mathbb{E}\big(f\big|\{\emptyset,\Omega\}\big) = \mathbb{E}f \quad \textup{a.s.}
\end{equation}

The following properties of conditional expectations are standard; see \cite[Section~4.1]{Du19:book}.
The operator $f\mapsto\mathbb{E}(f|\mathcal{G})$ is linear and monotone.
If $f\in\textup{L}^1(\Omega,\mathcal{F},\mathbb{P})$ and $g\in\textup{L}^\infty(\Omega,\mathcal{G},\mathbb{P})$, $\mathcal{G}\subseteq\mathcal{F}$, then
\begin{equation}\label{eq:ce3}
\mathbb{E}(fg|\mathcal{G}) = \mathbb{E}(f|\mathcal{G}) g \quad \textup{a.s.}
\end{equation}
Next, for $p\in[1,\infty)$ and $f\in\textup{L}^p(\Omega,\mathcal{F},\mathbb{P})$ we have
\begin{equation}\label{eq:ce4}
\big|\mathbb{E}(f|\mathcal{G})\big|^p \leq \mathbb{E}\big(|f|^p\big|\mathcal{G}\big) \quad \textup{a.s.} \end{equation}
In particular, the conditional expectation is a (not necessarily strict) contraction in the $\textup{L}^p$ norms.
Finally, if $\mathcal{G}$ and $\mathcal{H}$ are two $\sigma$-algebras such that $\mathcal{F}\supseteq\mathcal{G}\supseteq\mathcal{H}$, then
\begin{equation}\label{eq:ce5}
\mathbb{E}\big(\mathbb{E}(f|\mathcal{G})\big|\mathcal{H}\big) = \mathbb{E}(f|\mathcal{H})
= \mathbb{E}\big(\mathbb{E}(f|\mathcal{H})\big|\mathcal{G}\big) \quad \textup{a.s.}
\end{equation}
for any integrable function $f$.

Now suppose that we are given a \emph{filtration} $(\mathcal{G}_m)_{m=0}^{\infty}$ of the probability space $(\Omega,\mathcal{F},\mathbb{P})$, i.e., a sequence of $\sigma$-algebras satisfying $\mathcal{G}_0\subseteq\mathcal{G}_1\subseteq\mathcal{G}_2\subseteq\cdots\subseteq\mathcal{F}$.
For shortness, let us denote the conditional expectation operator $f\mapsto\mathbb{E}(f|\mathcal{G}_m)$ simply by $\mathbb{E}_m$, for each index $m$.
Conditional expectations with respect to filtrations are widely studied in the literature on martingales; see for instance \cite[Chapter~4]{Du19:book}. We will need the following slightly nonstandard inequality in Subsections~\ref{subsec:simplexstructured} and \ref{subsec:structrees}.
We claim that for any nonnegative bounded measurable function $f$ and for nonnegative integers $m,m_1,m_2,\ldots,m_n$ satisfying $m\leq\min\{m_1,\ldots,m_n\}$ we have
\begin{equation}\label{eq:mineq1}
\mathbb{E}_m \big( f (\mathbb{E}_{m_1}f) (\mathbb{E}_{m_2}f) \cdots (\mathbb{E}_{m_n}f) \big) \geq (\mathbb{E}_m f)^{n+1} \quad \textup{a.s.}
\end{equation}
For the proof of \eqref{eq:mineq1} we can assume, without loss of generality, that $m\leq m_1\leq \cdots \leq m_n$. We use the mathematical induction on $k=0,1,\ldots,n-1$ to show
\begin{equation}\label{eq:mineqa}
\mathbb{E}_m \bigg( \Big( \prod_{i=1}^{k}\mathbb{E}_{m_i}f \Big) (\mathbb{E}_{m_{k+1}}f)^{n+1-k} \bigg) \geq (\mathbb{E}_m f)^{n+1} \quad \textup{a.s.}
\end{equation}
The case $k=n-1$ of \eqref{eq:mineqa} is precisely \eqref{eq:mineq1}, because we can use \eqref{eq:ce5} and \eqref{eq:ce3}, respectively, to equate their left hand sides:
\[ \mathbb{E}_m \bigg( f \Big( \prod_{i=1}^{n}\mathbb{E}_{m_i}f \Big) \bigg)
= \mathbb{E}_m \mathbb{E}_{m_n} \bigg( f \Big( \prod_{i=1}^{n}\mathbb{E}_{m_i}f \Big) \bigg)
= \mathbb{E}_m \Big( \Big( \prod_{i=1}^{n-1}\mathbb{E}_{m_i}f \Big) (\mathbb{E}_{m_n}f)^2 \Big) \quad \textup{a.s.} \]
The induction basis $k=0$ of \eqref{eq:mineqa} is just a consequence of \eqref{eq:ce4} and \eqref{eq:ce5}:
\[ \mathbb{E}_m \big( (\mathbb{E}_{m_1}f)^{n+1} \big)
\geq (\mathbb{E}_m \mathbb{E}_{m_1}f)^{n+1}
= (\mathbb{E}_m f)^{n+1} \quad \textup{a.s.} \]
For the induction step we only need to rewrite and estimate the left hand side of \eqref{eq:mineqa} as
\begin{align*}
& \mathbb{E}_m \mathbb{E}_{m_{k}} \bigg( \Big( \prod_{i=1}^{k}\mathbb{E}_{m_i}f \Big) (\mathbb{E}_{m_{k+1}}f)^{n+1-k} \bigg)
= \mathbb{E}_m \bigg( \Big( \prod_{i=1}^{k}\mathbb{E}_{m_i}f \Big) \,\mathbb{E}_{m_{k}} \big( (\mathbb{E}_{m_{k+1}}f)^{n+1-k} \big) \bigg) \\
& \geq \mathbb{E}_m \bigg( \Big( \prod_{i=1}^{k}\mathbb{E}_{m_i}f \Big) (\mathbb{E}_{m_{k}}\mathbb{E}_{m_{k+1}}f)^{n+1-k} \bigg)
= \mathbb{E}_m \bigg( \Big( \prod_{i=1}^{k-1}\mathbb{E}_{m_i}f \Big) (\mathbb{E}_{m_{k}}f)^{n+2-k} \bigg) \quad \textup{a.s.},
\end{align*}
where we used properties \eqref{eq:ce5}, \eqref{eq:ce3}, \eqref{eq:ce4}, and \eqref{eq:ce5} again, in that order.
This completes the inductive proof of \eqref{eq:mineqa} and thus also confirms \eqref{eq:mineq1}.

An immediate consequence of \eqref{eq:mineq1} combined with \eqref{eq:ce2}, \eqref{eq:ce5}, and \eqref{eq:ce4} is a scalar inequality
\begin{equation}\label{eq:mineq2}
\mathbb{E} \big( f (\mathbb{E}_{m_1}f) (\mathbb{E}_{m_2}f) \cdots (\mathbb{E}_{m_n}f) \big) \geq (\mathbb{E} f)^{n+1}
\end{equation}
for any nonnegative bounded measurable $f$ and nonnegative integers $m_1,m_2,\ldots,m_n$.
On the other hand, an easy generalization of \eqref{eq:mineq1} is
\begin{equation}\label{eq:mineq3}
\mathbb{E}_m \Big( f \prod_{i=1}^{n}\mathbb{E}_{m_i}f \Big) \geq \Big( \overbrace{ \prod_{\substack{i\\ m_i<m}} \mathbb{E}_{m_i}f }^{N\text{ factors}} \Big) (\mathbb{E}_m f)^{n+1-N} \quad \textup{a.s.},
\end{equation}
where $m$ is now a completely arbitrary nonnegative integer and $N$ is the number of indices $i$ satisfying $m_i<m$.
One only needs to use \eqref{eq:ce3} to factor out $N$ terms from the left hand side of \eqref{eq:mineq3} and then apply \eqref{eq:mineq1} to the remaining $n+1-N$ terms.
Inequality \eqref{eq:mineq3} will be used to resolve nested conditional expectations when bounding them from below.

\section{Anisotropic simplices: proof of Theorem~\ref{thm:simplices}}
\label{sec:simplices}
Our approach is the closest in spirit to Bourgain's original proof from \cite{B86:roth}. One superficial difference is that we are using the heat kernel, where Bourgain used the Poisson kernel. Other notable differences are the treatment of ``mismatched'' scales in the structured part using \eqref{eq:Bourgainlower} below and the way in which nonlinear scales are treated in the error part. Even though Lyall and Magyar gave several very slick alternative proofs for the case of isotropic ``linear'' dilates of a nondegenerate simplex \cite{LM16:prod,LM19:distgraphs,LM19:hypergraphs} (also see \cite{HLM16:dist}), we still find Bourgain's proof the easiest one to adapt to our general scheme.

Recall that we were given positive numbers \eqref{eq:asandbs} and unit vectors \eqref{eq:justus}. However, here we embed $\mathbb{R}^n\cong\mathbb{R}^n\times\{0\}\subset\mathbb{R}^{n+1}$, so that the $u_k$ are now viewed as unit vectors in $\mathbb{R}^{n+1}$.
Let $\mu$ be the normalized Haar measure on the special orthogonal group $\textup{SO}(n+1,\mathbb{R})$.
For a compactly supported measurable function $f\colon\mathbb{R}^{n+1}\to[0,1]$ the counting form is defined as
\[ \mathcal{N}^{0}_{\lambda}(f) :=
\int_{\mathbb{R}^{n+1}} \int_{\textup{SO}(n+1,\mathbb{R})} f(x) \Big( \prod_{k=1}^{n} f(x + \lambda^{a_k} b_k U u_k) \Big) \,\textup{d}\mu(U) \,\textup{d}x, \]
while its smoothened variant is
\[ \mathcal{N}^{\varepsilon}_{\lambda}(f) :=
\int_{\mathbb{R}^{n+1}} \int_{\textup{SO}(n+1,\mathbb{R})} f(x) \Big( \prod_{k=1}^{n} (f\ast\mathbbm{g}_{(\varepsilon\lambda)^{a_k}b_k})(x + \lambda^{a_k} b_k U u_k) \Big) \,\textup{d}\mu(U) \,\textup{d}x \]
for $\lambda>0$ and $0<\varepsilon\leq 1$.
Denote
\begin{equation}\label{eq:defofaandc}
a:=a_1+a_2+\cdots+a_n,\quad c:=\min\{a_1,a_2,\ldots,a_n\}.
\end{equation}
As explained in Subsection~\ref{subsec:approach}, it is sufficient to show that for any choice of numbers $\delta,\varepsilon\in(0,1]$, positive integer $J$, scales $0<\lambda_1<\cdots<\lambda_J$ satisfying $\lambda_{j+1}\geq2\lambda_j$ for each index $j$, yet another scale $\lambda\in(0,\lambda_J]$, a sufficiently large number $R>0$ (depending on $J$ and the scales $\lambda_j$), and a measurable set $B\subseteq[0,R]^{n+1}$ with measure $|B|\geq \delta R^{n+1}$ we have
\begin{align}
& \mathcal{N}^{1}_{\lambda}(\mathbbm{1}_B) \gtrsim \delta^{n+1} R^{n+1}, \label{eq:simplexstructured} \\
& \sum_{j=1}^{J} \big|\mathcal{N}^{\varepsilon}_{\lambda_j}(\mathbbm{1}_B)-\mathcal{N}^{1}_{\lambda_j}(\mathbbm{1}_B)\big| \lesssim \varepsilon^{-(n+2)a} \log(1/\varepsilon) J^{1/2} R^{n+1}, \label{eq:simplexerror} \\
& \big|\mathcal{N}^{0}_{\lambda}(\mathbbm{1}_B)-\mathcal{N}^{\varepsilon}_{\lambda}(\mathbbm{1}_B)\big| \lesssim \varepsilon^{c/2} R^{n+1}. \label{eq:simplexuniform}
\end{align}
Saying that $R$ is sufficiently large means, for instance, assuming that $R\geq 2\lambda_J^{a_k}b_k$ for each $k\in\{1,\ldots,n\}$.

Once we have \eqref{eq:simplexstructured}--\eqref{eq:simplexuniform}, the actual argument establishing Theorem~\ref{thm:simplices} takes $\delta$ to be a fixed positive number smaller than $\overline{\delta}(A)$. Afterwards we choose $\varepsilon$ small enough such that
\[ \varepsilon^{c/2} \leq \vartheta_1 \delta^{n+1} \]
and then $J$ large enough so that
\[ \varepsilon^{-(n+2)a} \log(1/\varepsilon) J^{-1/2} \leq \vartheta_2 \delta^{n+1}. \]
Here $\vartheta_1,\vartheta_2>0$ are sufficiently small constants, depending only on the implicit constants from \eqref{eq:simplexstructured}--\eqref{eq:simplexuniform}.
Next, we take an unbounded sequence of scales $(\lambda_j)_{j=1}^{\infty}$ such that the set $A$ does not contain the desired configuration (which is here an anisotropic simplex) of size $\lambda_j$ for any index $j$.
As a consequence, $\mathcal{N}^{0}_{\lambda_j}(\mathbbm{1}_B)$ vanishes for each $j$.
We sparsify the sequence to achieve $\lambda_{j+1}\geq2\lambda_j$.
Take $x\in\mathbb{R}^{n+1}$ and a sufficiently large $R$ for which the set $B := (A-x)\cap[0,R]^{n+1}$ has measure at least $\delta R^{n+1}$.
By pigeonholing over the summands in \eqref{eq:simplexerror}, we can find an index $j\in\{1,2,\ldots,J\}$ such that $\mathcal{N}^{0}_{\lambda_j}(\mathbbm{1}_B)$ is at least a positive constant multiple of $\delta^{n+1} R^{n+1}$.
That way we arrive at a contradiction with $\mathcal{N}^{0}_{\lambda_j}(\mathbbm{1}_B)=0$.

\subsection{The structured part: proof of \eqref{eq:simplexstructured}}
\label{subsec:simplexstructured}
For $k=2,\ldots,n$ let us write the orthogonal projection of $u_k$ onto $\mathop{\textup{span}}(\{u_1,\ldots,u_{k-1}\})$ explicitly as
\[ \beta_{k,1} u_1 + \beta_{k,2} u_2 + \cdots + \beta_{k,k-1} u_{k-1} \]
for some scalars $\beta_{k,1},\beta_{k,2},\ldots,\beta_{k,k-1}\in\mathbb{R}$. Moreover, denote
\[ d_k := \mathop{\textup{dist}}\big(u_k, \mathop{\textup{span}}(\{u_1,\ldots,u_{k-1}\})\big) >0. \]
By the symmetry present in $\mathcal{N}^{\varepsilon}_{\lambda}$ when integrating over all systems $(y_1,\ldots,y_n)=$\linebreak $(\lambda^{a_1} b_1 U u_1,\ldots,\lambda^{a_n} b_n U u_n)$ we can rewrite the smoother counting form as in \cite{B86:roth}:
\begin{align*}
\mathcal{N}^{\varepsilon}_{\lambda}(f) = \int_{(\mathbb{R}^{n+1})^{n+1}}
& f(x) \,\Big( \prod_{k=1}^{n} (f\ast\mathbbm{g}_{(\varepsilon\lambda)^{a_k}b_k})(x + y_k) \Big)
\,\textup{d}\sigma^{y_1,\ldots,y_{n-1}}_{\lambda^{a_n}b_n}(y_n) \\
& \textup{d}\sigma^{y_1,\ldots,y_{n-2}}_{\lambda^{a_{n-1}}b_{n-1}}(y_{n-1}) \cdots
\textup{d}\sigma^{y_1}_{\lambda^{a_{2}}b_{2}}(y_2)
\,\textup{d}\sigma_{\lambda^{a_1}b_1}(y_1)
\,\textup{d}x.
\end{align*}
Here $\sigma$ denotes the spherical measure supported on $\mathbb{S}^n\subset\mathbb{R}^{n+1}$, while, for $k=2,\ldots,n$, we write $\sigma^{y_1,\ldots,y_{k-1}}$ for the $(n-k+1)$-dimensional spherical measure supported on the sphere that is centered at
\[ \beta_{k,1} \frac{y_1}{|y_1|} + \beta_{k,2} \frac{y_2}{|y_2|} + \cdots + \beta_{k,k-1} \frac{y_{k-1}}{|y_{k-1}|}, \]
has radius $d_k$, and belongs to an $(n-k+2)$-dimensional plane in $\mathbb{R}^{n+1}$ orthogonal to $\mathop{\textup{span}}(\{y_1,\ldots,y_{k-1}\})$.
We normalize these measures in a way that each of them has its total mass equal to $1$ and we view them as being defined on all Borel subsets of $\mathbb{R}^{n+1}$, even though their supports are lower-dimensional sets contained in the standard unit sphere $\mathbb{S}^n$.
Also note that the $n+1$ integrals over $\mathbb{R}^{n+1}$ are nested and the integration is performed in the order from the innermost one to the outermost one.

Here we concentrate on the smoothest case $\varepsilon=1$.
It is easy to observe that for any probability measure $\nu$ supported inside the standard unit sphere $\mathbb{S}^n\subset\mathbb{R}^{n+1}$ we can bound pointwise
\[ \nu\ast\mathbbm{g} \geq \Big(\inf_{\textup{B}(0,2)}\mathbbm{g}\Big) \mathbbm{1}_{\textup{B}(0,1)} \gtrsim \varphi, \]
where
\[ \varphi := |\textup{B}(0,1)|^{-1}\mathbbm{1}_{\textup{B}(0,1)}. \]
The integral in $y_n$ in $\mathcal{N}^{1}_{\lambda}(f)$ can be recognized as a triple convolution,
\[ \Big( f \ast \sigma^{y_1,\ldots,y_{n-1}}_{-\lambda^{a_n}b_n} \ast\mathbbm{g}_{\lambda^{a_n}b_n} \Big)(x)
\gtrsim (f\ast\varphi_{\lambda^{a_n}b_n})(x). \]
Now we do the same to the integral in $y_{n-1}$, etc. Repeating this process $n$ times we end up with
\[ \mathcal{N}^{1}_{\lambda}(f) \gtrsim \int_{(\mathbb{R}^{d})^{n+1}} f(x) \,\Big( \prod_{k=1}^{n} (f\ast\varphi_{\lambda^{a_k}b_k})(x) \Big) \,\textup{d}x. \]
Once we know that
\begin{equation}\label{eq:Bourgainlower}
\fint_{[0,R]^d} f(x) \,\Big( \prod_{k=1}^{n} (f\ast\varphi_{t_k})(x) \Big) \,\textup{d}x
\gtrsim \Big( \fint_{[0,R]^d} f(x) \,\textup{d}x \Big)^{n+1}
\end{equation}
holds for every $t_1,t_2,\ldots,t_n\in(0,R/2]$ and every measurable function $f\colon[0,R]^{d}\to[0,1]$, then \eqref{eq:simplexstructured} will follow simply by choosing $t_k=\lambda^{a_k}b_k$ for $k=1,\ldots,n$, and taking $d=n+1$ and $f=\mathbbm{1}_B$.
However, \eqref{eq:Bourgainlower} is just a straightforward generalization of \cite[Lemma~2.1]{DGR19:Roth} by Durcik, Guo, and Roos, which, in turn, is an elaboration of Bourgain's \cite[Lemma~6]{B88:nonlinroth} (stated there without proof). The paper \cite{DGR19:Roth} established the special case $d=1$, $n=2$, $R=1$ of \eqref{eq:Bourgainlower} elegantly, by dominating convolutions $f\ast\varphi_{t_k}$ pointwise from below by a dyadic martingale.
We are about to reuse this idea to give a quick proof of \eqref{eq:Bourgainlower} in general.
Even if this generalization is quite clear and expected, we still choose to be sufficiently detailed, since we will have to argue similarly in an even greater generality in Subsection~\ref{subsec:structrees}.

Let us turn $[0,R]^d$ into a probability space by taking $\mathbb{P}$ to be the Lebesgue measure normalized by the factor $R^{-d}$.
The \emph{dyadic filtration} $(\mathcal{G}_m)_{m=0}^{\infty}$ of $[0,R]^d$ is obtained by choosing $\mathcal{G}_m$ to be a finite $\sigma$-algebra (i.e., a finite algebra of sets) generated with $2^{dm}$ congruent cubes of sidelength $2^{-m}R$ that partition $[0,R]^d$.
Let $m_k$ be the smallest nonnegative integer such that $2^{-m_k}R d^{1/2} < t_k$.
By this choice a cube of sidelength $2^{-m_k}R$ has diameter smaller than $t_k$, so, if it contains a point $x$, then it remains fully inside the ball $B(x,t_k)$.
Because of this and Equation \eqref{eq:ce1} we clearly have
\[ f\ast\varphi_{t_k} \gtrsim \mathbb{E}(f|\mathcal{G}_{m_k}) = \mathbb{E}_{m_k}f \quad \textup{a.e.} \]
Now \eqref{eq:Bourgainlower} becomes a consequence of the probabilistic inequality \eqref{eq:mineq2}.

\subsection{The error part: proof of \eqref{eq:simplexerror}}
\label{subsec:simplexerror}
We will keep writing $f$ interchangeably with $\mathbbm{1}_B$, where $B$ is as before.
Using the product rule for differentiation and applying the generalized heat equation \eqref{eq:heateqgen}, we get
\begin{align*}
& \frac{\partial}{\partial t} \prod_{k=1}^{n} (f\ast\mathbbm{g}_{t^{a_k}\lambda^{a_k}b_k})(x_k) \\
& = \sum_{m=1}^{n} \frac{a_m}{2\pi t} (f\ast\mathbbm{k}_{t^{a_m}\lambda^{a_m}b_m})(x_m)
\Big( \prod_{\substack{1\leq k\leq n\\ k\neq m}} (f\ast\mathbbm{g}_{t^{a_k}\lambda^{a_k}b_k})(x_k) \Big).
\end{align*}
Thus, by the fundamental theorem of calculus, for any $0<\alpha<\beta$ the difference $\mathcal{N}^{\alpha}_{\lambda}(f)-\mathcal{N}^{\beta}_{\lambda}(f)$ can be written as
\[ \sum_{m=1}^{n}\mathcal{L}_{\lambda}^{\alpha,\beta,m}(f), \]
where
\begin{align}
\mathcal{L}_{\lambda}^{\alpha,\beta,m}(f) := -\frac{a_m}{2\pi} & \int_{\alpha}^{\beta} \int_{\mathbb{R}^{n+1}} \int_{\textup{SO}(n+1,\mathbb{R})}
f(x) \,(f\ast\mathbbm{k}_{(t\lambda)^{a_m}b_m})(x + \lambda^{a_m} b_m U u_m) \nonumber \\
& \times \Big(  \prod_{\substack{1\leq k\leq n\\ k\neq m}} (f\ast\mathbbm{g}_{(t\lambda)^{a_k}b_k})(x + \lambda^{a_k} b_k U u_k) \Big) \,\textup{d}\mu(U) \,\textup{d}x \,\frac{\textup{d}t}{t} . \label{eq:Bourgainformrep0}
\end{align}
By symmetry it is sufficient to consider the case $m=n$ and the same reasoning as in the previous subsection gives
\begin{align}
\mathcal{L}_{\lambda}^{\alpha,\beta,n}(f) = -\frac{a_n}{2\pi} \int_{\alpha}^{\beta} \int_{(\mathbb{R}^{n+1})^{n+1}}
& f(x) \,\Big( \prod_{k=1}^{n-1} (f\ast\mathbbm{g}_{(t\lambda)^{a_k}b_k})(x + y_k) \Big)  \nonumber \\
& \times (f\ast\mathbbm{k}_{(t\lambda)^{a_n}b_n})(x + y_n) \,\textup{d}\sigma^{y_1,\ldots,y_{n-1}}_{\lambda^{a_n}b_n}(y_n) \nonumber \\
& \cdots \textup{d}\sigma^{y_1}_{\lambda^{a_{2}}b_{2}}(y_2) \,\textup{d}\sigma_{\lambda^{a_1}b_1}(y_1)
\,\textup{d}x \,\frac{\textup{d}t}{t}. \label{eq:Bourgainformrep}
\end{align}

The following arguments will not depend on the dimension of the ambient space $\mathbb{R}^{n+1}$. We will emphasize this fact by writing it as $\mathbb{R}^d$ and remembering that $d=n+1\geq 2$. This will also be convenient for recycling the same computation in Section~\ref{sec:trees}.
Here we need to control
\[ \sum_{j=1}^{J}|\mathcal{L}_{\lambda_j}^{\varepsilon,1,n}(f)|. \]
Set
\begin{equation}\label{eq:defoftheta}
\theta:=10^{-1/a_n}e^{-1}.
\end{equation}
Let us begin the estimation by multiplying the inner integral of $\mathcal{L}_{\lambda_j}^{\varepsilon,1,n}(f)$ with
\begin{equation}\label{eq:integralis1}
\int_{\theta t\lambda_j}^{e\theta t\lambda_j} \,\frac{\textup{d}s}{s} = 1
\end{equation}
and rewriting the whole expression as
\begin{align*}
-\frac{a_n}{2\pi}\int_{\varepsilon}^{1} \int_{\theta t\lambda_j}^{e\theta t\lambda_j} \int_{(\mathbb{R}^{d})^{n}}
& f(x) \,\Big( \prod_{k=1}^{n-1} (f\ast\mathbbm{g}_{(t\lambda_j)^{a_k}b_k})(x + y_k) \Big) \\
& \times \big(f\ast\sigma^{y_1,\ldots,y_{n-1}}_{-\lambda_j^{a_n}b_n}\ast\mathbbm{k}_{(t\lambda_j)^{a_n}b_n}\big)(x) \,\textup{d}\sigma^{y_1,\ldots,y_{n-2}}_{\lambda_j^{a_{n-1}}b_{n-1}}(y_{n-1}) \\
& \cdots \textup{d}\sigma^{y_1}_{\lambda_j^{a_{2}}b_{2}}(y_2) \,\textup{d}\sigma_{\lambda_j^{a_1}b_1}(y_1)
\,\textup{d}x \,\frac{\textup{d}s}{s} \,\frac{\textup{d}t}{t}.
\end{align*}
For
\begin{equation}\label{eq:jstcond}
j\in\{1,\ldots,J\},\quad \varepsilon\leq t\leq 1,\quad \theta t\lambda_j\leq s\leq e\theta t\lambda_j
\end{equation}
denote
\[ r = r(j,s,t) := \big((t\lambda_j)^{2a_n}-s^{2a_n}\big)^{1/2a_n}, \]
remembering that $r$ depends on $j$, $s$, and $t$, and observing that $s\sim t\lambda_j \sim r$.
Convolution identity \eqref{eq:cvid2} gives
\[ \sigma^{y_1,\ldots,y_{n-1}}_{-\lambda_j^{a_n}}\ast\mathbbm{k}_{(t\lambda_j)^{a_n}}
= \Big(\frac{(t\lambda_j)^{2}}{rs}\Big)^{a_n} \sum_{l=1}^{d} \sigma^{y_1,\ldots,y_{n-1}}_{-\lambda_j^{a_n}} \ast \mathbbm{h}_{r^{a_n}}^{(l)} \ast \mathbbm{h}_{s^{a_n}}^{(l)}, \]
so that, introducing the integration variable $y$,
\begin{align*}
|\mathcal{L}_{\lambda_j}^{\varepsilon,1,n}(f)|
\lesssim \sum_{l=1}^{d} \int_{\varepsilon}^{1} & \int_{\theta t\lambda_j}^{e\theta t\lambda_j} \int_{(\mathbb{R}^{d})^{n+1}}
f(x) \,\Big( \prod_{k=1}^{n-1} (f\ast\mathbbm{g}_{(t\lambda_j)^{a_k}b_k})(x + y_k) \Big) \\
& \times \big|\big(f\ast\mathbbm{h}^{(l)}_{s^{a_n}b_n}\big)(x+y)\big| \,\big| \big(\sigma^{y_1,\ldots,y_{n-1}}_{\lambda_j^{a_n}b_n}\ast\mathbbm{h}^{(l)}_{r^{a_n}b_n}\big)(y) \big| \\
& \textup{d}\sigma^{y_1,\ldots,y_{n-2}}_{\lambda_j^{a_{n-1}}b_{n-1}}(y_{n-1}) \cdots \textup{d}\sigma_{\lambda_j^{a_1}b_1}(y_1)
\,\textup{d}x \,\textup{d}y \,\frac{\textup{d}s}{s} \,\frac{\textup{d}t}{t} .
\end{align*}
Take $\nu$ to be an arbitrary probability measure supported on a subset of the standard unit sphere in $\mathbb{R}^{d}$, even though the following reasoning will only be applied with $\nu=\sigma^{y_1,\ldots,y_{n-1}}$ in the current section.
Expanding the convolution according to its definition and using \eqref{eq:GaussianbySchwartz}, \eqref{eq:SchwartzbyGaussian} enables us to estimate
\begin{align*}
& \big|(\nu_{(\lambda_j/s)^{a_n}}\ast\mathbbm{h}^{(l)}_{(r/s)^{a_n}})(x)\big| \\
& \leq \int_{\mathbb{R}^{d}} \Big(\frac{s}{r}\Big)^{d a_n} \Big|\mathbbm{h}^{(l)}\Big(\Big(\frac{s}{r}\Big)^{a_n}x-\Big(\frac{\lambda_j}{r}\Big)^{a_n}y\Big)\Big| \,\textup{d}\nu(y) \\
& \lesssim \int_{\mathbb{R}^{d}} \Big(\frac{s}{r}\Big)^{d a_n} \Big(1+
\Big(\frac{s}{r}\Big)^{a_n}|x|\Big)^{-d-1} \Big(1+\Big(\frac{\lambda_j}{r}\Big)^{a_n}|y|\Big)^{d+1} \,\textup{d}\nu(y) \\
& \lesssim \varepsilon^{-(d+1) a_n} (1+|x|)^{-d-1}
\lesssim \varepsilon^{-(d+1) a_n} \int_1^\infty \mathbbm{g}_{\gamma}(x) \,\frac{\textup{d}\gamma}{\gamma^2},
\end{align*}
so, dilating by $s^{a_n}b_n$, we also get
\begin{equation}\label{eq:Gaussdomh}
\big|\big(\nu_{\lambda_j^{a_n}b_n}\ast\mathbbm{h}^{(l)}_{r^{a_n}b_n}\big)(x)\big| \lesssim \varepsilon^{-(d+1) a_n}\int_1^\infty \mathbbm{g}_{s^{a_n}b_n\gamma}(x) \,\frac{\textup{d}\gamma}{\gamma^2}
\end{equation}
for $j,t,s$ as in \eqref{eq:jstcond} and for $x\in\mathbb{R}^d$.
Consequently,
\begin{align*}
|\mathcal{L}_{\lambda_j}^{\varepsilon,1,n}(f)|
\lesssim \varepsilon^{-(d+1) a_n} & \sum_{l=1}^{d} \int_{1}^{\infty} \int_{\varepsilon}^{1} \int_{\theta t\lambda_j}^{e\theta t\lambda_j}  \int_{(\mathbb{R}^{d})^{n+1}}
f(x) \,\mathbbm{g}_{s^{a_n}b_n\gamma}(y) \\
& \times \big|\big(f\ast\mathbbm{h}^{(l)}_{s^{a_n}b_n}\big)(x+y)\big| \,\Big( \prod_{k=1}^{n-1} (f\ast\mathbbm{g}_{(t\lambda_j)^{a_k}b_k})(x + y_k) \Big) \\
& \textup{d}\sigma^{y_1,\ldots,y_{n-2}}_{\lambda_j^{a_{n-1}}b_{n-1}}(y_{n-1}) \cdots \textup{d}\sigma_{\lambda_j^{a_1}b_1}(y_1)
\,\textup{d}x \,\textup{d}y \,\frac{\textup{d}s}{s} \,\frac{\textup{d}t}{t} \,\frac{\textup{d}\gamma}{\gamma^2}.
\end{align*}

In the next step we realize that $y_{n-1}$ only appears in one of the factors of the integrand and as one of the integration variables in the last expression. For that reason, we can rewrite the integral in $y_{n-1}$ as a convolution and get
\begin{align}
|\mathcal{L}_{\lambda_j}^{\varepsilon,1,n}(f)|
\lesssim & \ \varepsilon^{-(d+1) a_n} \sum_{l=1}^{d} \int_{1}^{\infty} \int_{\varepsilon}^{1} \int_{\theta t\lambda_j}^{e\theta t\lambda_j}  \int_{(\mathbb{R}^{d})^{n}}
f(x) \,\mathbbm{g}_{s^{a_n}b_n\gamma}(y) \nonumber \\
& \times \big|\big(f\ast\mathbbm{h}^{(l)}_{s^{a_n}b_n}\big)(x+y)\big| \,\Big( \prod_{k=1}^{n-2} (f\ast\mathbbm{g}_{(t\lambda_j)^{a_k}b_k})(x + y_k) \Big) \nonumber \\
& \times \big(f\ast\sigma^{y_1,\ldots,y_{n-2}}_{-\lambda_j^{a_{n-1}}b_{n-1}}\ast\mathbbm{g}_{(t\lambda_j)^{a_{n-1}}b_{n-1}}\big)(x) \nonumber \\
& \,\textup{d}\sigma^{y_1,\ldots,y_{n-3}}_{\lambda_j^{a_{n-2}}b_{n-2}}(y_{n-2}) \cdots  \,\textup{d}\sigma_{\lambda_j^{a_1}b_1}(y_1)
\,\textup{d}x \,\textup{d}y \,\frac{\textup{d}s}{s} \,\frac{\textup{d}t}{t} \,\frac{\textup{d}\gamma}{\gamma^2}. \label{eq:splxformmed}
\end{align}
Now we observe that, by \eqref{eq:GaussianbySchwartz} and \eqref{eq:SchwartzbyGaussian},
\begin{align*}
& (\nu_{(\lambda_j/s)^{a_k}}\ast\mathbbm{g}_{(t\lambda_j/s)^{a_k}})(x) \\
& = \int_{\mathbb{R}^d} \Big(\frac{s}{t\lambda_j}\Big)^{d a_k} \mathbbm{g}\Big(\Big(\frac{s}{t\lambda_j}\Big)^{a_k}x-\frac{1}{t^{a_k}}y\Big) \,\textup{d}\nu(y) \\
& \lesssim \int_{\mathbb{R}^d} \Big(\frac{s}{t\lambda_j}\Big)^{d a_k} \Big(1+
\Big(\frac{s}{t\lambda_j}\Big)^{a_k}|x|\Big)^{-d-1} \Big(1+\frac{|y|}{t^{a_k}}\Big)^{d+1} \,\textup{d}\nu(y) \\
& \lesssim \varepsilon^{-(d+1) a_k} (1+|x|)^{-d-1}
\lesssim \varepsilon^{-(d+1) a_k} \int_1^\infty \mathbbm{g}_{\eta}(x) \,\frac{\textup{d}\eta}{\eta^2}
\end{align*}
for a general $\nu$ as before and for any $k\in\{1,\ldots,n-1\}$.
Rescaling this by $s^{a_k}b_k$ yields
\begin{equation}\label{eq:Gaussdomg}
(\nu_{\lambda_j^{a_k}b_k}\ast\mathbbm{g}_{(t\lambda_j)^{a_k}b_k})(x)
\lesssim \varepsilon^{-(d+1) a_k}\int_1^\infty \mathbbm{g}_{s^{a_k}b_k\eta}(x) \,\frac{\textup{d}\eta}{\eta^2}.
\end{equation}
Convolutions with dilates of $\mathbbm{g}$ can be controlled pointwise by the Hardy-Littlewood maximal function $\textup{M}$ (see \cite[Formula~(3.9)]{SW71:book}), so \eqref{eq:Gaussdomg} implies
\begin{equation}\label{eq:Gaussdomgmax}
(f\ast\nu_{\lambda_j^{a_k}b_k}\ast\mathbbm{g}_{(t\lambda_j)^{a_k}b_k})(x)
\lesssim \varepsilon^{-(d+1) a_k} (\textup{M}f)(x).
\end{equation}
We use \eqref{eq:Gaussdomgmax} with $k=n-1$ and $\nu=\sigma^{y_1,\ldots,y_{n-2}}_{-1}$ to further bound the expression \eqref{eq:splxformmed}.
Repeating the previous step $n-2$ times more, i.e., for $k=n-2,\ldots,2,1$, we end up with
\begin{align*}
|\mathcal{L}_{\lambda_j}^{\varepsilon,1,n}(f)|
\lesssim \varepsilon^{-(d+1) a} \sum_{l=1}^{d} & \int_{1}^{\infty} \int_{\varepsilon}^{1} \int_{\theta t\lambda_j}^{e\theta t\lambda_j}  \int_{(\mathbb{R}^{d})^{2}}
f(x) \,\big( (\textup{M}f)(x) \big)^{n-1} \\
& \times \big|\big(f\ast\mathbbm{h}^{(l)}_{s^{a_n}b_n}\big)(x+y)\big| \,\mathbbm{g}_{s^{a_n}b_n\gamma}(y)
\,\textup{d}x \,\textup{d}y \,\frac{\textup{d}s}{s} \,\frac{\textup{d}t}{t} \,\frac{\textup{d}\gamma}{\gamma^2} ,
\end{align*}
where $a$ is the sum of the numbers $a_k$, as in \eqref{eq:defofaandc}.
Observing $f\leq\textup{M}f$ and using the Cauchy--Schwarz inequality, we get
\begin{align*}
|\mathcal{L}_{\lambda_j}^{\varepsilon,1,n}(f)|^2
\lesssim & \ \varepsilon^{-2(d+1) a} \big(\log(1/\varepsilon)\big)
\|\textup{M}f\|_{\textup{L}^{2n}}^{2n} \\
& \times \Big( \int_{\varepsilon}^{1} \int_{\theta t\lambda_j}^{e\theta t\lambda_j}
\sum_{l=1}^{d} \big\|f\ast\mathbbm{h}^{(l)}_{s^{a_n}b_n}\big\|_{\textup{L}^2}^2
\,\frac{\textup{d}s}{s} \,\frac{\textup{d}t}{t} \Big).
\end{align*}
Now we sum in $j$ and observe that, for each fixed $t\in[\varepsilon,1]$, the intervals $[\theta t\lambda_j,e\theta t\lambda_j]$, $j=1,\ldots,J$, cover any fixed point from $(0,\infty)$ at most two times.
By one last application of the Cauchy--Schwarz inequality (this time for the sum in $j$) followed by boundedness of $\textup{M}$ on $\textup{L}^{2n}(\mathbb{R}^d)$, we conclude
\begin{align*}
\sum_{j=1}^{J} |\mathcal{L}_{\lambda_j}^{\varepsilon,1,n}(f)|
& \leq J^{1/2} \Big( \sum_{j=1}^{J} |\mathcal{L}_{\lambda_j}^{\varepsilon,1,n}(f)|^2 \Big)^{1/2} \\
& \lesssim J^{1/2} \varepsilon^{-(d+1) a} \big(\log(1/\varepsilon)\big)
\|f\|_{\textup{L}^{2n}}^{n}
\Big( \int_{0}^{\infty} \sum_{l=1}^{d} \big\|f\ast\mathbbm{h}^{(l)}_{s^{a_n}b_n}\big\|_{\textup{L}^2}^2 \,\frac{\textup{d}s}{s} \Big)^{1/2}.
\end{align*}
That way, identity \eqref{eq:idt1} and trivial estimates $\|f\|_{\textup{L}^{2n}}\leq R^{d/2n}$, $\|f\|_{\textup{L}^{2}}\leq R^{d/2}$, together with $d=n+1$, finish the proof of \eqref{eq:simplexerror}.

\subsection{The uniform part: proof of \eqref{eq:simplexuniform}}
\label{subsec:simplexuniform}
Take $0<\vartheta<\varepsilon$ and a measurable function $f\colon[0,R]^{n+1}\to[0,1]$. From the previous subsection we know that $\mathcal{N}^{\vartheta}_{\lambda}(f)-\mathcal{N}^{\varepsilon}_{\lambda}(f)$ is the sum of $\mathcal{L}_{\lambda}^{\vartheta,\varepsilon,m}(f)$ over $m=1,\ldots,n$.
Once again, by symmetry it is sufficient to consider the case $m=n$, when we have the representation \eqref{eq:Bourgainformrep}.
Using $0\leq f\leq 1$ and the Cauchy--Schwarz inequality in the variable $x$, we get
\begin{align*}
\big|\mathcal{L}_{\lambda}^{\vartheta,\varepsilon,n}(f)\big| \lesssim \int_{\vartheta}^{\varepsilon} \int_{(\mathbb{R}^{n+1})^{n-1}}
& \|f\|_{\textup{L}^2}
\big\| f \ast \sigma^{y_1,\ldots,y_{n-1}}_{-\lambda^{a_n}b_n}\ast\mathbbm{k}_{(t\lambda)^{a_n}b_n} \big\|_{\textup{L}^2} \\
& \textup{d}\sigma^{y_1,\ldots,y_{n-2}}_{\lambda^{a_{n-1}}b_{n-1}}(y_{n-1})
\cdots \textup{d}\sigma^{y_1}_{\lambda^{a_{2}}b_{2}}(y_2)
\,\textup{d}\sigma_{\lambda^{a_1}b_1}(y_1) \,\frac{\textup{d}t}{t}.
\end{align*}
Another application of the Cauchy--Schwarz inequality, this time in the variables $y_1,\ldots,$ $y_{n-1}$, followed by Plancherel's identity, leads to
\begin{align*}
\big|\mathcal{L}_{\lambda}^{\vartheta,\varepsilon,n}(f)\big| \lesssim  \|f\|_{\textup{L}^2(\mathbb{R}^{n+1})} \int_{\vartheta}^{\varepsilon} \Big( & \int_{(\mathbb{R}^{n+1})^{n}}
\big|\widehat{f}(\xi)\big|^2  \big|\widehat{\sigma}^{y_1,\ldots,y_{n-1}}(-\lambda^{a_n}b_n\xi)\big|^2 \\
& \times \big|\widehat{\mathbbm{k}}(t^{a_n}\lambda^{a_n}b_n\xi)\big|^2 \,\textup{d}\xi \,\textup{d}\sigma^{y_1,\ldots,y_{n-2}}_{\lambda^{a_{n-1}}b_{n-1}}(y_{n-1}) \\
& \cdots \textup{d}\sigma^{y_1}_{\lambda^{a_{2}}b_{2}}(y_2)
\,\textup{d}\sigma_{\lambda^{a_1}b_1}(y_1) \Big)^{1/2} \,\frac{\textup{d}t}{t}.
\end{align*}
Observe that $\sigma^{y_1,\ldots,y_{n-1}}$ is just the circle measure inside a two-dimensional plane orthogonal to $\mathop{\textup{span}}(\{y_1,\ldots,y_{n-1}\})$ as long as $y_1,\ldots,y_{n-1}$ are linearly independent, which happens almost surely.
Using the decay estimate \eqref{eq:subcircledecay} we get
\[ \big|\widehat{\sigma}^{y_1,\ldots,y_{n-1}}(\zeta)\big| \lesssim \mathop{\textup{dist}}\big(\zeta,\mathop{\textup{span}}(\{y_1,\ldots,y_{n-1}\})\big)^{-1/2} \]
for every $\zeta\in\mathbb{R}^{n+1}$. By rescaling $y_1,\ldots,y_{n-1}$ this in turn gives
\begin{align}
& \big|\mathcal{L}_{\lambda}^{\vartheta,\varepsilon,n}(f)\big| \nonumber \\
& \lesssim \|f\|_{\textup{L}^2(\mathbb{R}^{n+1})} \int_{\vartheta}^{\varepsilon} \Big( \int_{\mathbb{R}^{n+1}}
\big|\widehat{f}(\xi)\big|^2 \big|\widehat{\mathbbm{k}}(t^{a_n}\lambda^{a_n}b_n\xi)\big|^2
\mathcal{I}(\lambda^{a_n}b_n\xi) \,\textup{d}\xi \Big)^{1/2} \,\frac{\textup{d}t}{t}, \label{eq:Bourgainunif1}
\end{align}
where
\begin{align*}
\mathcal{I}(\zeta) := \int_{(\mathbb{R}^{n+1})^{n-1}} & \mathop{\textup{dist}}\big(\zeta,\mathop{\textup{span}}(\{y_1,\ldots,y_{n-1}\})\big)^{-1} \\
& \textup{d}\sigma^{y_1,\ldots,y_{n-2}}(y_{n-1}) \cdots \textup{d}\sigma^{y_1}(y_2) \,\textup{d}\sigma(y_1).
\end{align*}

Integrating over all rotations $U\in \textup{SO}(n+1,\mathbb{R})$ and taking $(y_1,\ldots,y_{n-1})=(U u_1,\ldots,$ $U u_{n-1})$ we want to conclude
\begin{equation}\label{eq:estfortermI}
\mathcal{I}(\zeta) \lesssim |\zeta|^{-1}.
\end{equation}
We could proceed as in \cite{B86:roth} or \cite{LM16:prod}, but a more elementary argument is also available.
Instead of integrating over $U$, we can rather take $y_1,\ldots,y_{n-1}$ to be fixed unit vectors that span the coordinate plane $\mathbb{R}^{n-1}\times\{(0,0)\}$ and integrate over all possible directions of $\zeta$ determined by the standard unit sphere $\mathbb{S}^{n}$ in $\mathbb{R}^{n+1}$. This will yield the same result for $\mathcal{I}(\zeta)$ up to a constant. Writing $\zeta$ in the $(n+1)$-dimensional spherical coordinates,
\[ \zeta = |\zeta| (\cos\phi_1, \sin\phi_1\cos\phi_2, \ldots, \sin\phi_1\cdots\sin\phi_{n-1}\cos\phi_n, \sin\phi_1\cdots\sin\phi_n), \]
$\phi_1,\ldots,\phi_{n-1}\in[0,\pi]$, $\phi_n\in[0,2\pi)$, and observing
\[ \mathop{\textup{dist}}\big(\zeta,\mathop{\textup{span}}(\{y_1,\ldots,y_{n-1}\})\big) = |\zeta| \sin\phi_1 \cdots \sin\phi_{n-1}, \]
we can estimate $\mathcal{I}(\zeta)$ by a constant multiple of
\[ |\zeta|^{-1} \int_{[0,\pi]^{n-1}\times[0,2\pi)} \frac{\sin^{n-1}\phi_1\sin^{n-2}\phi_2 \cdots \sin\phi_{n-1} \,\textup{d}\phi_1\,\textup{d}\phi_2\cdots\textup{d}\phi_n}{\sin\phi_1\sin\phi_2 \cdots \sin\phi_{n-1}} \lesssim |\zeta|^{-1}. \]
This establishes \eqref{eq:estfortermI} and thus also gives
\begin{equation}\label{eq:Bourgainunif2}
|\widehat{\mathbbm{k}}(t^{a_n}\zeta)|^2 \mathcal{I}(\zeta) \lesssim (t^{a_n}|\zeta|)^4 e^{-2\pi (t^{a_n}|\zeta|)^2} |\zeta|^{-1} \lesssim t^{a_n}
\end{equation}
for $\zeta\in\mathbb{R}^{n+1}$ and $t>0$.
Substituting $\zeta=\lambda^{a_n}b_n\xi$, plugging \eqref{eq:Bourgainunif2} into \eqref{eq:Bourgainunif1}, and using Plancherel's theorem again, we finally get
\[ \big|\mathcal{L}_{\lambda}^{\vartheta,\varepsilon,n}(f)\big|
\lesssim \|f\|_{\textup{L}^2(\mathbb{R}^{n+1})}^2 \int_{\vartheta}^{\varepsilon} t^{a_n/2} \,\frac{\textup{d}t}{t}
\lesssim R^{n+1}\varepsilon^{a_n/2}. \]
This proves \eqref{eq:simplexuniform}.

Arguments in this subsection reveal irrelevance of the nature of dilations for this part of the proof, which is the reason why we were able to proceed in a similar way as Bourgain \cite{B86:roth} or Lyall and Magyar \cite{LM19:distgraphs,LM16:prod}.

\section{Anisotropic rectangular boxes: proof of Theorem~\ref{thm:boxes}}
\label{sec:boxes}
In this section $\sigma$ will denote exclusively the circular measure in $\mathbb{R}^2$. Many elements of the proof will be similar to the corresponding ingredients in \cite[Section~7]{DK20:Szemeredi}. Still, a few modifications are needed.

For a compactly supported measurable function $f\colon(\mathbb{R}^2)^n\to[0,1]$ and for $\lambda>0$ and $0<\varepsilon\leq 1$ this time we define
\[ \mathcal{N}^{0}_{\lambda}(f) :=
\int_{(\mathbb{R}^2)^{2n}} \Big( \prod_{(r_1,\ldots,r_n)\in\{0,1\}^n} \!\!\! f(x_1 + r_1 y_1, \ldots, x_n + r_n y_n) \Big)
\Big( \prod_{k=1}^{n}\textup{d}x_k\,\textup{d}\sigma_{\lambda^{a_k}b_k}(y_k) \Big) \]
and
\begin{align*}
\mathcal{N}^{\varepsilon}_{\lambda}(f)
& := \int_{(\mathbb{R}^2)^{2n}} \Big( \prod_{(r_1,\ldots,r_n)\in\{0,1\}^n} f(x_1 + r_1 y_1, \ldots, x_n + r_n y_n) \Big) \\
& \qquad\qquad\ \times \Big( \prod_{k=1}^{n}(\sigma\ast\mathbbm{g}_{\varepsilon^{a_k}})_{\lambda^{a_k}b_k}(y_k)\,\textup{d}x_k\,\textup{d}y_k \Big) \\
& \,= \int_{(\mathbb{R}^2)^{2n}} \mathcal{F}(\mathbf{x})\,
\Big( \prod_{k=1}^{n} (\sigma\ast\mathbbm{g}_{\varepsilon^{a_k}})_{\lambda^{a_k}b_k}(x_k^0-x_k^1) \Big) \,\textup{d}\mathbf{x},
\end{align*}
where we recall the notation \eqref{eq:boxauxdef1}--\eqref{eq:boxauxdef2}, specialized to $d=2$.

Theorem~\ref{thm:boxes} will follow once we establish:
\begin{align}
& \mathcal{N}^{1}_{\lambda}(\mathbbm{1}_B) \gtrsim \delta^{2^n} R^{2n}, \label{eq:cubesstructured} \\
& \sum_{j=1}^{J} \big|\mathcal{N}^{\varepsilon}_{\lambda_j}(\mathbbm{1}_B)-\mathcal{N}^{1}_{\lambda_j}(\mathbbm{1}_B)\big| \lesssim \varepsilon^{-3a} \log(1/\varepsilon) R^{2n}, \label{eq:cubeserror} \\
& \big|\mathcal{N}^{0}_{\lambda}(\mathbbm{1}_B)-\mathcal{N}^{\varepsilon}_{\lambda}(\mathbbm{1}_B)\big| \lesssim \varepsilon^{c/2} R^{2n}. \label{eq:cubesuniform}
\end{align}
Here, $a$, $c$, $\delta$, $\varepsilon$, $J$, $\lambda_j$, $\lambda$, and $R$ are just as in Section~\ref{sec:simplices}, while $B\subseteq([0,R]^2)^n$ is a measurable set satisfying $|B|\geq\delta R^{2n}$.
We set $f=\mathbbm{1}_B$ and continue using the notation \eqref{eq:boxauxdef1}--\eqref{eq:boxauxdef2}.

\subsection{The structured part: proof of \eqref{eq:cubesstructured}}
\label{subsec:boxesstructured}
We partition a major part of the cube $([0,R]^2)^n$ into the collection of rectangular boxes $Q_1 \times \cdots \times Q_n$, where each $Q_k$ is a square of the form $[l \lambda^{a_k}b_k,(l+1) \lambda^{a_k}b_k)\times [l' \lambda^{a_k}b_k,(l'+1) \lambda^{a_k}b_k)$ for some integers $0\leq l,l'\leq \lfloor\lambda^{-a_k}b_k^{-1}R\rfloor-1$.
Each of these boxes has measure $(\lambda^{a_k}b_k)^2$ and their total number is clearly comparable to $R^{2n}\prod_{k=1}^{n}\lambda^{-2a_k}$.

When estimating $\mathcal{N}^{1}_{\lambda}(\mathbbm{1}_B)$ from below, we restrict the domain of integration with additional constraints requiring that $x_k^0$, $x_k^1$ lie in the same square $Q_k$ mentioned above, for each $k=1,\ldots,n$. Using the box--Gowers--Cauchy--Schwarz inequality \cite{S06:corners,GT08:primes,Tao07:exp}, or simply by several applications of the ordinary Cauchy--Schwarz inequality, we obtain
\begin{align*}
\fint_{Q_1\times Q_1\times\cdots\times Q_n\times Q_n}
\mathcal{F}(\mathbf{x})\,\textup{d}\mathbf{x}
\geq \Big(\fint_{Q_1\times\ldots\times Q_n} f\Big)^{2^n}.
\end{align*}
Applying this to each of the choices of $Q_1\times\cdots\times Q_n$, using
\[ (\sigma\ast\mathbbm{g})_{\lambda^{a_k}b_k} \gtrsim \lambda^{-2a_k} \mathbbm{1}_{[-\lambda^{a_k}b_k,\lambda^{a_k}b_k]^2}, \]
and recalling $|Q_k|\sim\lambda^{2a_k}$, we get
\[  \mathcal{N}^{1}_{\lambda}(\mathbbm{1}_B) \gtrsim \Big(\prod_{k=1}^{n}\lambda^{-2a_k}\Big) \Big(\prod_{k=1}^{n}\lambda^{4a_k}\Big) \sum_{Q_1\times\cdots\times Q_n} \Big(\fint_{Q_1\times\ldots\times Q_n} \mathbbm{1}_B\Big)^{2^n}, \]
so discrete Jensen's inequality for the power function $t\mapsto t^{2^n}$ gives \eqref{eq:cubesstructured}.

\subsection{The error part: proof of \eqref{eq:cubeserror}}
\label{subsec:boxeserror}
Using the product rule and the generalized heat equation \eqref{eq:heateqgen} we obtain
\begin{align*}
\frac{\partial}{\partial t} \prod_{k=1}^{n} (\sigma_{\lambda^{a_k}b_k}\ast\mathbbm{g}_{t^{a_k}\lambda^{a_k}b_k})(y_k)
= \sum_{m=1}^{n} \frac{a_m}{2\pi t}  & (\sigma_{\lambda^{a_m}b_m}\ast\mathbbm{k}_{t^{a_m}\lambda^{a_m}b_m})(y_m) \\
& \times \Big( \prod_{\substack{1\leq k\leq n\\ k\neq m}} (\sigma_{\lambda^{a_k}b_k}\ast\mathbbm{g}_{t^{a_k}\lambda^{a_k}b_k})(y_k) \Big)
\end{align*}
for $\lambda>0$ and $y_1,\ldots,y_n\in\mathbb{R}^2$.
By the fundamental theorem of calculus the difference $\mathcal{N}^{\alpha}_{\lambda}(f)-\mathcal{N}^{\beta}_{\lambda}(f)$ is, for any $0<\alpha<\beta$, equal to the sum of $n$ terms given by
\begin{align*}
\mathcal{L}_{\lambda}^{\alpha,\beta,m}(f) := -\frac{a_m}{2\pi} \int_{\alpha}^{\beta} \int_{(\mathbb{R}^2)^{2n}}
& \mathcal{F}(\mathbf{x}) \, (\sigma\ast\mathbbm{k}_{t^{a_m}})_{\lambda^{a_m}b_m}(x_m^0-x_m^1) \\
& \times \Big( \prod_{\substack{1\leq k\leq n\\ k\neq m}} (\sigma\ast\mathbbm{g}_{t^{a_k}})_{\lambda^{a_k}b_k}(x_k^0-x_k^1) \Big) \,\textup{d}\mathbf{x} \,\frac{\textup{d}t}{t}
\end{align*}
for $m=1,\ldots,n$. By symmetry it is sufficient to fix $m=n$ and prove an upper bound for
\begin{equation}\label{eq:cubeserrlhs}
\sum_{j=1}^{J}|\mathcal{L}_{\lambda_j}^{\varepsilon,1,n}(f)|.
\end{equation}

Take $\theta$ as in \eqref{eq:defoftheta} and use \eqref{eq:integralis1}, transforming $\mathcal{L}_{\lambda_j}^{\varepsilon,1,n}(f)$ into
\begin{align*}
-\frac{a_n}{2\pi}\int_{\varepsilon}^{1} \int_{\theta t\lambda_j}^{e\theta t\lambda_j} \int_{(\mathbb{R}^2)^{2n}}
& \mathcal{F}(\mathbf{x}) \,\Big( \prod_{k=1}^{n-1} (\sigma\ast\mathbbm{g}_{t^{a_k}})_{\lambda_j^{a_k}b_k}(x_k^0-x_k^1) \Big) \\
& \times (\sigma\ast\mathbbm{k}_{t^{a_n}})_{\lambda_j^{a_n}b_n}(x_n^0-x_n^1)
\,\textup{d}\mathbf{x} \,\frac{\textup{d}s}{s} \,\frac{\textup{d}t}{t}.
\end{align*}
For $j,s,t$ as in \eqref{eq:jstcond} this time we denote
\[ r = r(j,s,t) := \big((t\lambda_j)^{2a_n}-2s^{2a_n}\big)^{1/2a_n} \]
and observe, once again, that $s\sim t\lambda_j \sim r$.
From \eqref{eq:cvid3} and \eqref{eq:cvid2} we see
\[ \sigma_{\lambda_j^{a_n}}\ast\mathbbm{k}_{(t\lambda_j)^{a_n}}
= \Big(\frac{t\lambda_j}{s}\Big)^{2a_n} \sum_{l=1}^2 \sigma_{\lambda_j^{a_n}} \ast \mathbbm{g}_{r^{a_n}} \ast \mathbbm{h}_{s^{a_n}}^{(l)} \ast \mathbbm{h}_{s^{a_n}}^{(l)}, \]
so \eqref{eq:cubeserrlhs} is at most a constant times
\begin{align*}
\sum_{j=1}^{J} \sum_{l=1}^2 \int_{\varepsilon}^{1} \int_{\theta t\lambda_j}^{e\theta t\lambda_j} \int_{(\mathbb{R}^2)^{2n}}
& \Big| \int_{\mathbb{R}^2} \mathcal{F}^{(n)}(x_n^0) \,\mathbbm{h}_{s^{a_n}b_n}^{(l)}(x_n^0-q^0) \,\textup{d}x_n^0 \Big| \\
& \times \Big| \int_{\mathbb{R}^2} \mathcal{F}^{(n)}(x_n^1) \,\mathbbm{h}_{s^{a_n}b_n}^{(l)}(x_n^1-q^1) \,\textup{d}x_n^1 \Big| \\
& \times \Big( \prod_{k=1}^{n-1} (\sigma_{\lambda_j^{a_k}b_k} \ast \mathbbm{g}_{(t\lambda_j)^{a_k}b_k})(x_k^0-x_k^1) \,\textup{d}x_{k}^0 \,\textup{d}x_{k}^1 \Big) \\
& \times (\sigma_{\lambda_j^{a_n}b_n} \ast \mathbbm{g}_{r^{a_n}b_n})(q^0-q^1) \,\textup{d}q^0 \,\textup{d}q^1 \,\frac{\textup{d}s}{s} \,\frac{\textup{d}t}{t} .
\end{align*}
The same computation leading to \eqref{eq:Gaussdomg} still applies, so we can use this formula again, this time with $d=2$ and $\nu=\sigma$:
\begin{equation}\label{eq:Gaussdomgrep1}
(\sigma_{\lambda_j^{a_k}b_k}\ast\mathbbm{g}_{(t\lambda_j)^{a_k}b_k})(x)
\lesssim \varepsilon^{-3 a_k}\int_1^\infty \mathbbm{g}_{s^{a_k}b_k\gamma}(x) \,\frac{\textup{d}\gamma}{\gamma^2}
\end{equation}
for $k=1,\ldots,n-1$.
Very similarly, imitating the proof of \eqref{eq:Gaussdomh}, we get
\begin{equation}\label{eq:Gaussdomgrep2}
(\sigma_{\lambda_j^{a_n}b_n}\ast\mathbbm{g}_{r^{a_n}b_n})(x) \lesssim \varepsilon^{-3a_n}\int_1^\infty \mathbbm{g}_{s^{a_n}b_n\gamma}(x) \,\frac{\textup{d}\gamma}{\gamma^2}.
\end{equation}
Estimates \eqref{eq:Gaussdomgrep1} and \eqref{eq:Gaussdomgrep2} bound \eqref{eq:cubeserrlhs} by a constant multiple of
\begin{align*}
\varepsilon^{-3a} & \sum_{j=1}^{J} \sum_{l=1}^2 \int_{[1,\infty)^n} \int_{\varepsilon}^{1} \int_{\theta t\lambda_j}^{e\theta t\lambda_j} \int_{(\mathbb{R}^2)^{2n}}
\Big| \int_{\mathbb{R}^2} \mathcal{F}^{(n)}(x_n^0) \,\mathbbm{h}_{s^{a_n}b_n}^{(l)}(x_n^0-q^0) \,\textup{d}x_n^0 \Big| \\
& \times \Big| \int_{\mathbb{R}^2} \mathcal{F}^{(n)}(x_n^1) \,\mathbbm{h}_{s^{a_n}b_n}^{(l)}(x_n^1-q^1) \,\textup{d}x_n^1 \Big|
\,\Big( \prod_{k=1}^{n-1} \mathbbm{g}_{s^{a_k} b_k \gamma_k}(x_k^0-x_k^1) \,\textup{d}x_k^0 \,\textup{d}x_k^1 \Big) \\
& \times \mathbbm{g}_{s^{a_n} b_n \gamma_n}(q^0-q^1) \,\textup{d}q^0 \,\textup{d}q^1 \,\frac{\textup{d}s}{s} \,\frac{\textup{d}t}{t} \,\frac{\textup{d}\gamma_1}{\gamma_1^2} \cdots \frac{\textup{d}\gamma_n}{\gamma_n^2} .
\end{align*}
Using the Cauchy--Schwarz inequality in all variables other than $x_n^0$ and $x_n^1$, observing that the two obtained terms are equal, expanding out the square, and collapsing back the convolution using identity \eqref{eq:cvid2}, we see that this expression is at most
\begin{align}
-\frac{1}{2} \varepsilon^{-3a} \int_{[1,\infty)^n} \int_{\varepsilon}^{1} \sum_{j=1}^{J} \int_{\theta t\lambda_j}^{e\theta t\lambda_j} \int_{(\mathbb{R}^2)^{2n}}
\mathcal{F}(\mathbf{x})
\,\Big( \prod_{k=1}^{n-1} \mathbbm{g}_{s^{a_k} b_k \gamma_k}(x_k^0-x_k^1) \Big) & \nonumber \\
\times \mathbbm{k}_{2^{1/2}s^{a_n}b_n}(x_n^0-x_n^1)
\,\textup{d}\mathbf{x} \,\frac{\textup{d}s}{s} \,\frac{\textup{d}t}{t} \,\frac{\textup{d}\gamma_1}{\gamma_1^2} \cdots \frac{\textup{d}\gamma_n}{\gamma_n^2} & . \label{eq:simplifedcubeexp}
\end{align}

For a fixed $t$ and varying $j$ the number of times the intervals $[\theta t\lambda_j,e\theta t\lambda_j]$ cover any fixed point from $(0,\infty)$ is at most a constant depending only on $a_n$.
Using this observation in connection with \eqref{eq:simplifedcubeexp} and recognizing the inner expression as \eqref{eq:defofthetaform}, we finally obtain
\[ \sum_{j=1}^{J}|\mathcal{L}_{\lambda_j}^{\varepsilon,1,n}(f)| \lesssim \varepsilon^{-3a} \log(1/\varepsilon) \int_{[1,\infty)^n}  \Theta^{n,n}_{\gamma_1,\ldots,\gamma_{n-1},2^{1/2}}(f) \,\frac{\textup{d}\gamma_1}{\gamma_1^2} \cdots \frac{\textup{d}\gamma_n}{\gamma_n^2}. \]
Thus, the bound for $\Theta^{n,n}_{\gamma_1,\ldots,\gamma_{n-1},2^{1/2}}(f)$ coming from identity \eqref{eq:infortheta} completes the proof of \eqref{eq:cubeserror}.

\subsection{The uniform part: proof of \eqref{eq:cubesuniform}}
\label{subsec:cubesuniform}
Take $0<\vartheta<\varepsilon\leq 1$. In order to control $\mathcal{N}^{\vartheta}_{\lambda}(f)-\mathcal{N}^{\varepsilon}_{\lambda}(f)$ we need to bound $|\mathcal{L}_{\lambda}^{\vartheta,\varepsilon,m}(f)|$ for $m=1,\ldots,n$.

Once we fix $m$ and a number $t$ such that $\vartheta\leq t\leq \varepsilon$, Plancherel's theorem yields
\begin{align}
& \int_{(\mathbb{R}^2)^{2}} \mathcal{F} \, (\sigma\ast\mathbbm{k}_{t^{a_m}})_{\lambda^{a_m}b_m}(x_m^0-x_m^1) \,\textup{d}x_m^0 \,\textup{d}x_m^1 \nonumber \\
& = \int_{\mathbb{R}^2} \big|\widehat{\mathcal{F}^{(m)}}(\xi)\big|^2 \,\widehat{\sigma}(\lambda^{a_m}b_m\xi) \,\widehat{\mathbbm{k}}(t^{a_m}\lambda^{a_m}b_m\xi) \,\textup{d}\xi. \label{eq:unifcubes}
\end{align}
We use the decay estimate \eqref{eq:circledecay} for $\widehat{\sigma}$ to conclude
\begin{equation}\label{eq:unifcubesdecay}
\big|\widehat{\sigma}(\zeta)\widehat{\mathbbm{k}}(t^{a_m}\zeta)\big|
\lesssim |\zeta|^{-1/2} (t^{a_m}|\zeta|)^2 e^{-\pi (t^{a_m}|\zeta|)^2} \lesssim t^{a_m/2}
\end{equation}
for each $\zeta\in\mathbb{R}^2$. Then we simply take $\zeta=\lambda^{a_m}b_m\xi$ in \eqref{eq:unifcubesdecay} and combine it with \eqref{eq:unifcubes}.
By another application of Plancherel's theorem we see that \eqref{eq:unifcubes} is at most a constant times
\[ t^{a_m/2} \|\mathcal{F}^{(m)}\|_{\textup{L}^2(\mathbb{R}^2)}^2 \leq t^{a_m/2} R^2. \]
Recall that $f=\mathbbm{1}_B$.
Integrating in all of the remaining variables we get
\[ | \mathcal{L}_{\lambda}^{\vartheta,\varepsilon,m}(\mathbbm{1}_B) | \lesssim \int_{\vartheta}^{\varepsilon} t^{a_m/2} R^{2n} \,\frac{\textup{d}t}{t} \]
and thus also
\[ \big|\mathcal{N}^{\vartheta}_{\lambda}(\mathbbm{1}_B)-\mathcal{N}^{\varepsilon}_{\lambda}(\mathbbm{1}_B)\big| \lesssim
\int_{\vartheta}^{\varepsilon} t^{c/2} R^{2n} \,\frac{\textup{d}t}{t}. \]
It remains to send $\vartheta\to0^+$.

\section{Anisotropic trees: proof of Theorem~\ref{thm:trees}}
\label{sec:trees}
After all of the material presented in Sections~\ref{sec:simplices} and \ref{sec:boxes} the proof of Theorem~\ref{thm:trees} will not require many new ideas. We merely need to pick and reapply a few elements of the proofs of Theorems~\ref{thm:simplices} and \ref{thm:boxes}, so we will be somewhat brief.

Once again, let $\sigma$ be the circular measure supported on $\mathbb{S}^1\subseteq\mathbb{R}^2$. Take a compactly supported measurable function $f\colon\mathbb{R}^2\to[0,1]$ and parameters $\lambda\in(0,\infty)$, $\varepsilon\in(0,1]$. Let us jump straight to the definition of the smoother counting form:
\begin{align*}
\mathcal{N}^{\varepsilon}_{\lambda}(f)
& := \int_{(\mathbb{R}^2)^{n+1}}
\Big( \prod_{k\in E} (\sigma\ast\mathbbm{g}_{\varepsilon^{a_k}})_{\lambda^{a_k}b_k}(x_{u(k)}-x_{v(k)}) \Big)
\Big( \prod_{v\in V} f(x_v) \textup{d}x_v \Big),
\end{align*}
where we write $u(k),v(k)$ for the two vertices that are joined by the edge $k\in E$.
The actual counting form $\mathcal{N}^{0}_{\lambda}$, obtained in the limit as $\varepsilon\to0$, can also be defined explicitly by declaring a particular vertex as the ``root'' of the tree; see a similar reasoning in Subsection~\ref{subsec:structrees} below. We do not need a formula for $\mathcal{N}^{0}_{\lambda}(f)$ in the proof.

In order to prove Theorem~\ref{thm:trees} it is enough to show the following:
\begin{align}
& \mathcal{N}^{1}_{\lambda}(\mathbbm{1}_B) \gtrsim \delta^{n+1} R^{2}, \label{eq:treesstructured} \\
& \sum_{j=1}^{J} \big|\mathcal{N}^{\varepsilon}_{\lambda_j}(\mathbbm{1}_B)-\mathcal{N}^{1}_{\lambda_j}(\mathbbm{1}_B)\big| \lesssim \varepsilon^{-3a} \log(1/\varepsilon) J^{1/2} R^{2}, \label{eq:treeserror} \\
& \big|\mathcal{N}^{0}_{\lambda}(\mathbbm{1}_B)-\mathcal{N}^{\varepsilon}_{\lambda}(\mathbbm{1}_B)\big| \lesssim \varepsilon^{c/2} R^{2}, \label{eq:treesuniform}
\end{align}
where $a$, $c$, $J$, $\lambda_j$, $\lambda$, $R$ were described in Section~\ref{sec:simplices} and $B\subseteq[0,R]^2$ is a measurable set satisfying $|B|\geq\delta R^2$. We write $f$ interchangeably with $\mathbbm{1}_B$.

\subsection{The structured part: proof of \eqref{eq:treesstructured}}
\label{subsec:structrees}
Denote $t_k:=\lambda^{a_k}b_k$ for each $k\in E$. Just as in Subsection~\ref{subsec:simplexstructured}, turn $[0,R]^2$ into a probability space and define the dyadic filtration $(\mathcal{G}_m)_{m=0}^{\infty}$ on it.

Let us declare arbitrary particular vertex $v_{\mathcal{T}}\in V$ as the \emph{tree root} or the \emph{tree top}. We can imagine that the tree $\mathcal{T}$ is ``hanged'' upside down by holding it by the vertex $v_{\mathcal{T}}$. This naturally yields to relations \emph{is a child of} and \emph{is a parent of} on the set of vertices $V$. For any $v\in V$ let $\mathcal{T}_v=(V_v,E_v)$ be the subtree of $\mathcal{T}$ consisting of $v$ as its root and of all descendants of $v$ with respect to $\mathcal{T}$.

For any tree $\mathcal{T}=(V,E)$ with root $v_{\mathcal{T}}$ we define the following function of only one variable $x_{v_\mathcal{T}}$:
\[ (\mathcal{A}_{\mathcal{T}}f)(x_{v_\mathcal{T}}) :=
\int_{(\mathbb{R}^2)^{|E|}} \Big( \prod_{k\in E} (\sigma\ast\mathbbm{g})_{t_k}(x_{u(k)}-x_{v(k)}) \Big)
\Big( \prod_{v\in V} f(x_v) \Big) \Big( \prod_{\substack{v\in V\\ v\neq v_{\mathcal{T}}}} \textup{d}x_v \Big). \]
By induction on the tree structure we are going to show that for any such tree $\mathcal{T}$ there exist nonnegative integers $m_1,m_2,\ldots,m_{|E|}$ (depending also on the numbers $t_k$) such that the following pointwise inequality holds:
\begin{equation}\label{eq:indtre}
\mathcal{A}_{\mathcal{T}}f \geq f (\mathbb{E}_{m_1}f) (\mathbb{E}_{m_2}f) \cdots (\mathbb{E}_{m_{|E|}}f) \quad \textup{a.e.}
\end{equation}
for any nonnegative bounded measurable function $f$.
Indeed, the induction basis is the case when the tree $\mathcal{T}$ has only one vertex and no edges, and then \eqref{eq:indtre} reduces to a trivial inequality $f\geq f$.
For the induction step let $k_1,\ldots,k_s$ be all edges incident with the root $v_\mathcal{T}$ and let $v_1,\ldots,v_s$ respectively be the corresponding children of $v_\mathcal{T}$.
Clearly,
\[ (\mathcal{A}_{\mathcal{T}}f)(x_{v_\mathcal{T}})
= \int_{(\mathbb{R}^2)^{s}} f(x_{v_{\mathcal{T}}}) \Big( \prod_{i=1}^{s} \big( \mathcal{A}_{\mathcal{T}_{v_i}}f \big)(x_{v_i}) \,(\sigma\ast\mathbbm{g})_{t_{k_i}}(x_{v_{\mathcal{T}}}-x_{v_i}) \,\textup{d}x_{v_i} \Big), \]
i.e.,
\[ \mathcal{A}_{\mathcal{T}}f
= f \,\prod_{i=1}^{s} \big( (\mathcal{A}_{\mathcal{T}_{v_i}}f) \ast (\sigma\ast\mathbbm{g})_{t_{k_i}} \big) \quad \textup{a.e.} \]
For each $i\in\{1,\ldots,s\}$ let $l_i$ be the smallest nonnegative integer such that $2^{-{l_i}}R<t_{k_i}$. That way we get
\begin{equation}\label{eq:indtrf}
\mathcal{A}_{\mathcal{T}}f \gtrsim f \ \prod_{i=1}^{s} \mathbb{E}_{l_i} \mathcal{A}_{\mathcal{T}_{v_i}}f \quad \textup{a.e.}
\end{equation}
Let us apply the induction hypothesis of \eqref{eq:indtre} to each of the subtrees $\mathcal{T}_{v_i}$ and then resolve the nested conditional expectations using inequality \eqref{eq:mineq3}.
Plugging these into \eqref{eq:indtrf} and multiplying them for all $i$, we conclude that $\mathcal{A}_{\mathcal{T}}f$ also satisfies a lower bound of the form \eqref{eq:indtre} for some positive integers $m_1,\ldots,m_{|E|}$, which we do not need to make explicit. This finalizes the induction step and establishes the claim \eqref{eq:indtre}.

Combining \eqref{eq:indtre} with \eqref{eq:mineq2} and applying it to $f=\mathbbm{1}_B$ we finally conclude
\begin{align*}
\mathcal{N}^{1}_{\lambda}(f)
& = \int_{\mathbb{R}^2} \mathcal{A}_{\mathcal{T}}f
\geq R^2 \fint_{[0,R]^2} f (\mathbb{E}_{m_1}f) (\mathbb{E}_{m_2}f) \cdots (\mathbb{E}_{m_{|E|}}f) \\
& \geq R^2 \Big(\fint_{[0,R]^2} f \Big)^{|V|} \geq \delta^{|V|} R^2,
\end{align*}
which is precisely \eqref{eq:treesstructured}.

\subsection{The error part: proof of \eqref{eq:treeserror}}
\label{subsec:errtrees}
Generalized heat equation \eqref{eq:heateqgen} and the fundamental theorem of calculus allow us to rewrite the difference $\mathcal{N}^{\varepsilon}_{\lambda}(f)-\mathcal{N}^{1}_{\lambda}(f)$ as
\[ \sum_{m\in E}\mathcal{L}_{\lambda}^{\varepsilon,1,m}(f), \]
where, this time,
\begin{align*}
\mathcal{L}_{\lambda}^{\alpha,\beta,m}(f) := -\frac{a_m}{2\pi} & \int_{\alpha}^{\beta} (\sigma\ast\mathbbm{k}_{t^{a_m}})_{\lambda^{a_m}b_m}(x_{u(m)}-x_{v(m)}) \\
& \times\Big( \prod_{\substack{k\in E\\ k\neq m}} (\sigma\ast\mathbbm{g}_{t^{a_k}})_{\lambda^{a_k}b_k}(x_{u(k)}-x_{v(k)}) \Big)
\Big( \prod_{v\in V} f(x_v) \textup{d}x_v \Big) \,\frac{\textup{d}t}{t}.
\end{align*}
The proof of \eqref{eq:simplexerror} presented in Subsection~\ref{subsec:simplexerror} does not recognize any graph structure at all. Thus, the same proof carries over here, replacing measures $\sigma^{y_1,\ldots,y_{k-1}}$ (associated with varying spheres of different dimensions) always with the same circle measure $\sigma$.

\subsection{The uniform part: proof of \eqref{eq:treesuniform}}
The proof of \eqref{eq:simplexuniform} given in Subsection~\ref{subsec:simplexuniform} does not see any graph structure either. It keeps cancellation at only one crucial place, i.e., associated with only one chosen edge. For this reason we can proceed with only very minor modifications of the arguments from either Subsection~\ref{subsec:simplexuniform} or Subsection~\ref{subsec:cubesuniform}.

\section{Closing remarks}
\label{sec:closing}

\subsection{Multilinear anisotropic singular integrals}
\label{subsec:anisosingint}
As we have already said, a part of the motivation behind this work lies in encouraging connections between the techniques from multilinear harmonic analysis and the problems on combinatorics of the Euclidean space.
In this subsection we want to single out singular integrals that can naturally be associated with problems studied in Theorems~\ref{thm:simplices}--\ref{thm:trees}.
This correspondence should not be understood literally, but rather on the level of heuristics and methodology.
After all, the proofs of the above theorems did not use any estimates for the integral operators that will be mentioned here. In Sections~\ref{sec:simplices}--\ref{sec:trees} we preferred to use ad hoc shortcuts in the form of the Gaussian domination estimate \eqref{eq:SchwartzbyGaussian} and identities \eqref{eq:idt1} and \eqref{eq:infortheta}.
However, it would be a pity not to mention analytical counterparts of these problems, if for nothing else, then to point out where to look for the techniques for handling similar or more general combinatorial questions.

First, for a tuple $\mathbf{f}=(f_0,f_1,\ldots,f_n)$ of measurable functions on $\mathbb{R}^d$ and for a certain singular kernel $K$ we define the translation-invariant multilinear form
\begin{equation}\label{eq:opersimpl}
\Lambda_{K}(\mathbf{f}) := \textup{p.v.} \int_{(\mathbb{R}^{d})^{n+1}} K(x_1-x_0,\ldots,x_n-x_0) \,\Big(\prod_{k=0}^{n} f_k(x_k)\,\textup{d}x_k \Big).
\end{equation}
From either \eqref{eq:Bourgainformrep0} or \eqref{eq:Bourgainformrep} we are naturally lead to the study of multilinear singular integral forms \eqref{eq:opersimpl}.
For instance, if we were allowed to replace the measures $\mu$ or $\sigma^{y_1,\ldots,y_{k-1}}$ by the Dirac delta measure at the origin, we would obtain the kernel $K$ given by
\begin{equation}\label{eq:kersimpl}
K(y_1,\ldots,y_n) := \int_{0}^{\infty} \Big( \prod_{k=1}^{n-1} \mathbbm{g}_{t^{a_k}b_k}(y_k) \Big) \mathbbm{k}_{t^{a_n}b_n}(y_n) \,\frac{\textup{d}t}{t}.
\end{equation}
The study of multilinear singular integrals \eqref{eq:opersimpl} with Calder\'{o}n-Zygmund kernels $K$ was initiated by Coifman and Meyer in the 1970s (for instance see \cite{CM75:mCZ,CM78:mCZ1,CM78:mCZ2}), while a more systematic treatment was first given by Grafakos and Torres \cite{GT02:mCZ}.
However, \eqref{eq:kersimpl} is not the most usual singular kernel. It satisfies the Calder\'{o}n-Zygmund estimates with respect to quasinorms associated with anisotropic power-type dilations \eqref{eq:anidilsm}. Such more general dilation structures have been studied by Stein and Wainger \cite{SW78:dil}.
Strictly speaking, classical results, such as those from \cite{GT02:mCZ}, do not apply to kernels \eqref{eq:kersimpl}, but the same techniques do, and the proofs can be repeated mutatis mutandis (also see \cite{GBMS16:CZ}).
Unsurprisingly, Subsection~\ref{subsec:simplexerror} already comes quite close to proving some $\textup{L}^p$ bounds for \eqref{eq:opersimpl}, with its use of maximal and square function estimates, the trivial case of the latter being identity \eqref{eq:idt1}.
Hovewer, it is conceivable that, in the future, one reduces a certain combinatorial problem to the study of \eqref{eq:opersimpl} for quite different kernels $K$.

Next, for a tuple of measurable functions $\mathbf{f}=(f_{r_1,\ldots,r_n})_{(r_1,\ldots,r_n)\in\{0,1\}^n}$ on $(\mathbb{R}^d)^{n}$ indexed by points from $\{0,1\}^n$ and for a singular kernel $K$ we define
\begin{align}
\Theta_{K}(\mathbf{f}) := \textup{p.v.} \int_{(\mathbb{R}^d)^{2n}} & \Big( \prod_{(r_1,\ldots,r_n)\in\{0,1\}^n} f_{r_1,\ldots,r_n}(x_1 + r_1 y_1, \ldots, x_n + r_n y_n) \Big) \nonumber \\
& \times \,K(y_1,\ldots,y_n) \Big( \prod_{k=1}^{n}\textup{d}x_k\,\textup{d}y_k \Big) . \label{eq:operboxes}
\end{align}
These forms have also been studied extensively when $K$ is the usual Calder\'{o}n-Zygmund kernel, i.e., satisfying Calder\'{o}n-Zygmund estimates with respect to the Euclidean metric.
The first $\textup{L}^p$ bounds for the forms \eqref{eq:operboxes} were established in the case $d=1$ and $n=2$ by Durcik \cite{D15:L4,D17:Lp}. Prior to that, their dyadic model has been investigated by the present author \cite{K11:Bell,K12:twist} and by Thiele and the present author \cite{KT13:T1}. Estimates for rather general ``entangled'' multilinear singular integrals of the above type (i.e., with cubical structure) have been shown recently by Durcik and Thiele \cite{DT20:BL}.
It is also interesting to mention that the study of multiparameter variants of these objects has only recently been initiated by Bernicot and Durcik \cite{BD20:twist}.
It could be interesting to study \eqref{eq:operboxes} for more general dilation structures, i.e., when $K$ is a generalized Calder\'{o}n-Zygmund kernel, such as in the case \eqref{eq:kersimpl}, which is relevant here. In this context, a similar but still different object has been studied by \v{S}kreb and the present author \cite{KS15:mart}. Some generalizations of the result by Durcik \cite{D15:L4} are straightforward: the single $\textup{L}^{2^n}\times\cdots\times\textup{L}^{2^n}$ bound for \eqref{eq:operboxes} can be extracted easily from the presented proof of Theorem~\ref{thm:boxes} combined with a cone decomposition of the kernel.
It is quite likely that $\Theta_K$ still satisfies the same $\textup{L}^p$ estimates from \cite{D17:Lp}, but we do not attempt such generalizations here.

Finally, after Theorem~\ref{thm:trees}, a closely related topic of further investigation could be establishing estimates for entangled multilinear singular integrals associated with bipartite graphs or $r$-partite $r$-regular hypergraphs. Dyadic model of these problems are significantly easier: they have already been handled quite generally by the present author \cite{K11:Bell} and Stip\v{c}i\'{c} \cite{S20:T1}, respectively. Otherwise, the only cases and variants of entangled singular integral forms studied so far are the so-called ``twisted paraproduct operator'' \cite{K12:twist,DR18:averSHT}, the operators with cubical structure \cite{D15:L4,D17:Lp,DT20:BL}, and the operators that resemble multilinear Hilbert transforms \cite{T16:mht,Z17:splx,DKT16:splx,DK20:Szemeredi}.

\subsection{Comments on the anisotropic setting}
\label{subsec:comments}
Configurations generated by anisotropic power-type dilations \eqref{eq:anidilsm} have not been studied prior to this work. We found this setting sufficiently interesting because it fits nicely to the general method explained in Subsection~\ref{subsec:approach}.
Still, the results formulated in Section~\ref{subsec:newresults} are far from being definite.

In relation with simplices, the fact that we are dilating by power functions $\lambda\mapsto\lambda^a b$ plays little role in the proof of Theorem~\ref{thm:simplices}. Structured and uniform parts are handled much more generally, while the control of the error part essentially requires control of the multiscale objects of the form \eqref{eq:opersimpl}. Perhaps by studying multilinear singular integrals forms \eqref{eq:opersimpl} more closely one can hope for further generalizations of Theorem~\ref{thm:simplices}. In the present paper, the setting \eqref{eq:anidilsm} was just very convenient.

The same comment does not apply to rectangular boxes. So far, the only known way of handing entangled singular forms \eqref{eq:operboxes} is quite rigid and uses the same steps as those employed in the proof of Theorem~\ref{thm:boxes}. All of the papers \cite{D15:L4,D17:Lp,DKST19:NVEA,DKT16:splx,DT20:BL} used domination by Gaussians \eqref{eq:SchwartzbyGaussian} and some form of the product rule \eqref{eq:theprodrule}, or equivalently, integration by parts. This leaves less flexibility, so algebraic properties, such as \eqref{eq:cvid1}--\eqref{eq:heateqgen}, are now crucial in the proof.

Finally, one might wonder why we do not study general distance graphs in Theorem~\ref{thm:trees}. When vertices of the graph (or points of the configuration) are added one by one, it is desirable that the location of the new vertex depends on the locations of previous vertices and $\lambda$ in a reasonably simple way.
In general, dilations \eqref{eq:anidilsm} with different exponents $a_k$ can make this dependence ``very nonlinear.''
This complication arises already when the distance graph in question is a cycle of length $3$.
Suppose that, for all sufficiently large $\lambda$, we want to find points $x_0,x_1,x_2\in A\subseteq\mathbb{R}^d$ such that $|x_0-x_1|=\lambda$, $|x_0-x_2|=|x_1-x_2|=\lambda^2$.
Once $x_0$ and $x_1$ are located, the third point $x_2$ has to lie on the hyperplane $H$ orthogonal to the segment joining $x_0$ and $x_1$ and passing through its midpoint $(x_0+x_1)/2$. However, the distance from $x_3$ to $(x_0+x_1)/2$ needs to be $\lambda\sqrt{4\lambda^2-1}/2$ and we would need to stretch the corresponding spherical measure by this radical factor.
For a similar reason in Theorem~\ref{thm:simplices} we do not determine a simplex solely by lengths of its edges, but rather fix the angles between all pairs of its edges meeting at a single vertex.

\section*{Acknowledgments}
This work is supported in part by the \emph{Croatian Science Foundation} under project UIP-2017-05-4129 (MUNHANAP) and in part by the \emph{Fulbright Scholar Program}.
The author is grateful to Alex Amenta, Polona Durcik, and Jo\~{a}o Pedro Ramos for useful discussions and to both anonymous referees for numerous suggestions on improving the text. He also appreciates hospitality of the \emph{Georgia Institute of Technology} in the academic year 2019--20.


\end{document}